\documentclass[journal,twoside,web]{ieeecolor}
\usepackage{generic}
\usepackage{cite}
\usepackage{amsmath,amssymb,amsfonts}
\usepackage{algorithmic}
\usepackage{graphicx}
\usepackage{textcomp}

\usepackage{amsopn, amstext}
\usepackage{colortbl}
\usepackage{empheq}
\usepackage{graphicx}
\usepackage{color}
\usepackage[OT1]{fontenc}
\usepackage{mathtools}
\usepackage{scalerel}

\newtheorem{remark}{Remark}
\newtheorem{theorem}{Theorem}
\newtheorem{definition}{Definition}
\newtheorem{lemma}{Lemma}

\newtheorem{problem}{Problem}

\newtheorem{proposition}{Proposition}
\def\r{\color{black}}
\def\b{\color{black}}
\def\g{\color{black}}
\def\dg{\color{black}}
\def\rr{\color{black}}
\def\s{\small}
\def\n{\normalsize}
\def\f{\footnotesize}
\def\+{\!+\!}
\def\-{\!-\!}
\def\={\!=\!}
\def\cP{\mathcal{P}}
\def\cPm{\mathcal{P}_0}
\def\cPT{\mathcal{P}^{\!\top}\!}
\def\cPTm{\mathcal{P}_0^{\!\top}\!}
\def\E{\mathbb{E}}
\def\Eo{\E^{W_0}}

\def\z{\bar{z}}

\def\tu{\tilde{u}}
\def\F{\mathcal{F}}

\def\1mm{\vspace{-1mm}}
\def\2mm{\vspace{-2mm}}
\def\3mm{\vspace{-3mm}}
\def\4mm{\vspace{-4mm}}

\newenvironment{shrinkeq}[1]
{ \bgroup
\addtolength\abovedisplayshortskip{#1}
\addtolength\abovedisplayskip{#1}
\addtolength\belowdisplayshortskip{#1}
\addtolength\belowdisplayskip{#1}}
{\egroup\ignorespacesafterend}

\definecolor{deepgreen}{rgb}{0,0.5,0}
\def\BibTeX{{\rm B\kern-.05em{\sc i\kern-.025em b}\kern-.08em
    T\kern-.1667em\lower.7ex\hbox{E}\kern-.125emX}}
\begin{document}
\title{Mixed Social Optima and Nash equilibrium in Linear-Quadratic-Gaussian Mean-field System}
\author{Xinwei~Feng$^1$, Jianhui Huang$^2$ and Zhenghong~Qiu$^{3,*}$
  \thanks{$*$ Corresponding author}
  \thanks{$1$ \!Zhongtai Securities Institute for Financial Studies, Shandong University, Jinan, Shandong 250100, China. Email: xwfeng@sdu.edu.cn}
  \thanks{$2$ \!Department of Applied Mathematics, The Hong Kong Polytechnic University, Hong Kong. Email: majhuang@polyu.edu.hk}
  \thanks{$3$ \!Department of Applied Mathematics, The Hong Kong Polytechnic University, Hong Kong. Email: zhenghong.qiu@connect.polyu.hk}
}

\maketitle

\begin{abstract}
  This paper investigates a class of mixed stochastic
  linear-quadratic-Gaussian (LQG) social optimization and Nash game in
  the context of a large-scale system. Two types of interactive agents
  are involved: a major agent and a large number of weakly-coupled
  minor agents. All minor agents are cooperative to minimize the
  social cost as the sum of their individual costs, whereas such
  social cost is conflictive to that of the major agent. Thus, the
  major agent and all minor agents are further competitive to reach
  some  nonzero-sum Nash equilibrium.
   Applying the mean-field
   approximations and person-by-person optimality, we obtain auxiliary
   control problems for the major agent and minor agents,
   respectively. The decentralized social strategy is derived by a
   class of new  consistency condition (CC) system, which {\r consists of} mean-field
   forward-backward stochastic differential equations. The
   well-posedness of CC system is obtained by the discounting
   method. The related asymptotic social optimality for minor agents
   and Nash equilibrium for major-minor agents are also verified.
\end{abstract}


\section{Introduction}
{\r The Mean-field methodology} for large-population systems has been extensively
studied recently.
The central goal of an individual
 agent in the large-population system is to obtain decentralized strategy based on its own limited information.
  One efficient approach is the mean-field method which enables us to obtain
  the decentralized strategy through the limiting auxiliary control problem
  and the related consistency condition (CC) system. Along this direction, the
  interested readers are referred to  \cite{HCM2007}, \cite{HCM2015} and
  \cite{LL2007} for the derivation of mean-field games, \cite{BSYS2016},
  \cite{HHL2017}, \cite{HWW2016} for linear-quadratic-Gaussian (LQG)
  mean-field games, \cite{CD2013} for probabilistic analysis in mean-field
  games, \cite{basar2018} for risk-sensitive mean-field games,
  \cite{basar2015} for discrete-time mean-field games. \cite{NCMH2013} studied the mean-field social solution to consensus problems.  In the basic mean-field
  decision model, all agents have comparably small influence. However, in some
  real models, there exists an agent that has a significant influence on other
  agents. Thus a modified framework is to introduce a major agent interacting
  with a large number of minor agents. \cite{Huang2010} considered LQG major-minor games. \cite{NC2013}
  studies large population dynamic major-minor games involving nonlinear stochastic
  dynamical systems. 
   \cite{HHN2018} considered heterogeneous major-minor mean-field games. {\cite{carmona2017alternative} and \cite{carmona2016probabilistic}
  studied the probabilistic approach to major-minor mean-field games.}


In this paper, within the mean-field modeling, we will investigate  a new class of stochastic LQG optimization problems involving a major agent and a large number of weakly-coupled minor agents.
Specifically, the minor agents are cooperative to minimize the social cost as the sum of individual costs, while
the major agent and minor agents are competitive aiming for Nash equilibrium in nonzero-sum game manner.


   Our setup is an extension of the well-studied two-player (non-cooperative) game in which two agents make competitive decisions based on \emph{individual but centralized} information. As an extension, in our setup: one agent no longer applies such centralized decision. Instead, all its sub-units or branches will apply distributed information to  optimize the original cost jointly (e.g., \cite{B1996}, \cite{HK1983} and \cite{RV1999}), that now is reformulated in some team-cost form. Thereby, all sub-units become ``minor" agents and formalize a (cooperative) team, while another agent still applying centralized information becomes a ``major" and non-cooperative player, from the viewpoint of all ``minor" agents.


In our study, the problem can be solved in the following way. Firstly, for the major agent, we freeze the state average
 and obtain the auxiliary control problem. By the result in \cite{YZ1999}, the  auxiliary control problem for the major agent can be derived.
  Secondly, for the minor agents, under the person-by-person optimality principle, by applying variational techniques and introducing some mean-field terms, the original minor social optimization problem can also be converted to an auxiliary LQG control problem, which can be solved using some traditional scheme in \cite{YZ1999} as well. Thirdly, to determine the frozen mean-field terms, we construct the consistency condition (CC) system by some fixed-point analysis.
 Last, by using some asymptotic analysis and standard estimation of stochastic differential equations (SDEs), we show that the mean-field strategy provides an efficient approximation (i.e., the optimal loss tends to 0 as the population $N$ tends to $\infty$).

Moreover, the innovative aspects of the obtained results in this paper are as follows: Firstly, in general mean-field games framework, usually the auxiliary control problem can be obtained directly by replacing the state average with some frozen mean-field term. However,  this scheme will bring some ``bad" strategy in our social optima framework, which can not achieve the asymptotic optimality. Instead, in Section IV, variational techniques are applied to distinguish the high-order infinitesimals after the mean-filed approximation, and a new type of auxiliary control problem would be derived.
Secondly, the state process and state average enter the  diffusion terms. Such feature brings many difficulties when we apply the variational
method to obtain the auxiliary control problem for the minor agents. In particular, $N+1$ additional adjoint processes should be introduced to deal with the cross-terms in the cost functional variation.
 Thirdly, in the estimation of optimal loss, unlike the general mean-field games framework, the asymptotic optimality is proved through investigating the Fr\'{e}chet derivative of the social cost in Section VII.
Last, the control process enters the  diffusion terms. Because of this, the adjoint-state term will enter the drift term of the CC system. This also brings difficulties when we study the solvability of the CC system, which is a mean-field forward-backward stochastic differential equations (MF-FBSDEs) system, and the mean-field terms are represented in an embedding way. To its well-posedness, we apply some discounting methods.

The remaining of the paper is organized as follows:
In Section II, we give the formulation of the mixed LQG social optima problem. In Section III and Section IV, we
 find the auxiliary control problem of the major agent and minor agents respectively. The CC system is derived in Section V. Meanwhile, the well-posedness of  CC system is also established. In Section VI, we compare our result with some previous literature. In Section VII,  we obtain the asymptotic optimality of the decentralized strategy. Last, in Section VIII,  we simulate our model with some numerical methods.

\section{Problem formulation}\label{formulation}

Consider a finite time horizon $[0,T]$ for fixed $T>0$. Assume that $(\Omega,\F,\{\F_t\}_{0\leq t\leq T},\mathbb P)$ is a complete filtered probability space satisfying the usual conditions and $\{W_i(t),0\leq i\leq N\}_{0\leq t\leq T}$ is an $(N+1)$-dimensional Brownian motion on this space.  Let $\F_t$ be the filtration generated by $\{W_i(s),0\leq i\leq N\}_{0\leq s\leq t}$ and augmented by $\mathcal N_{\mathbb P}$ (the class of all $\mathbb P$-null sets of $\F$). Let $\F_t^{W_0}$, $\F_t^{W_i}$ and $\F_t^i$ be the augmentation of $\sigma\{W_0(s),0\leq s\leq t\}$, $\sigma\{W_i(s),0\leq s\leq t\}$ and $\sigma\{W_0(s),W_i(s),0\leq s\leq t\}$ by $\mathcal N_{\mathbb P}$ respectively. \rr $\Eo$ denotes the conditional expectation w.r.t. $\mathcal{F}^{W_0}_\cdot$. \b

Let $\langle\cdot,\cdot\rangle$ {\g denote} standard Euclidean inner product and $\|\cdot\|$ {\g denote} the norm. $x^\top$ denotes the transpose of a vector (or matrix) $x$. $\mathbb{S}^n$ denotes the set of symmetric $n\times n$ matrices with real elements. $M> (\geq) 0$ denotes that $M\in \mathbb{S}^n$ which is positive (semi)definite, while $M\gg 0$ denotes that, for some $\varepsilon>0$, $M - \varepsilon I \geq 0$. We introduce the following spaces for any given Euclidean space $\mathbb{H}$.  They will be used in the paper:
\begin{itemize} \s
  \item $L^2_{\F_T}(\Omega;\mathbb{H})\!:=\! \{\eta:\Omega\!\rightarrow\!\mathbb{H}|\eta$ is $\F_T$-measurable, $\E \|\eta\|^2\!<\!\infty\}$,
  \item $L_{\F_t}^2(0,T;\mathbb{H}):= \{\zeta(\cdot):[0,T]\times\Omega\rightarrow\mathbb{H}|\zeta(\cdot)$ is $\F_t$-progressively measurable,
        $\E \int_0^T\|\zeta({\g t})\|^2dt<\infty\}$.
  \item $L^\infty(0,\!T;\!\mathbb{H}) \!:= \! \{\zeta(\cdot)\!:\![0,\!T]\!\rightarrow\! \mathbb{H}|\text{esssup}_{0\leq s\leq T}\|\zeta(s)\|\!<\!\!\infty\}$.
  \item $L^2_{\F_t}(\Omega;C(0,T;\mathbb{H})) := \{\zeta(\cdot):[0,T]\times\Omega\rightarrow\mathbb{H}|\zeta(\cdot)$ is $\F_t$-adapted, continuous, $\E (\sup_{0\leq s\leq T} \|\zeta(s)\|^2)<\infty\}$,\n
\end{itemize}
and $\|\zeta\|^2_{L^2}:=\E \int_0^T\|\zeta\|^2dt$ denotes the $L^2$ norm. We consider a weakly coupled large population system with a major agent $\mathcal A_0$ and $N$ individual minor agents denoted by $\{\mathcal A_i:1\leq i\leq N\}$. The dynamics of the $N+1$ agents are given by a system of SDEs with mean-field coupling:
\begin{shrinkeq}{-1ex}
\begin{equation}\label{state equation-major}\resizebox{8.09cm}{!}{$\left\{\begin{aligned}
         & dx_0(t)=[A_0(t)x_0(t)+B_0(t)u_0(t)+F_0(t)x^{(N)}(t)]dt                         \\
         & \hspace{0.5cm}+[C_0(t)x_0(t)+D_0(t)u_0(t)+\widetilde F_0(t)x^{(N)}(t)]dW_0(t), \\
         & x_0(0)=\xi_0\in\mathbb R^n,                                                            \\
      \end{aligned}\right.$}\end{equation}\end{shrinkeq}
and for $1\leq i\leq N,$
\begin{shrinkeq}{-1ex}\begin{equation}\label{state equation-minor}
\resizebox{8.09cm}{!}{$  \left\{\begin{aligned}
     & dx_i(t)=[A(t)x_i(t)\+B(t)u_i(t)\+F(t)x^{(N)}(t)]dt                                                                           \\
     & \hspace{0.5cm}+\![C(t)x_i(t) \+D(t)u_i(t) \+ \widetilde  F(t)x^{(N)}\!(t) \+ \widetilde G(t)x_0(t) ] dW_i(t), \\ &x_i(0)=\xi\in\mathbb R^n,\\
  \end{aligned}\right.$} \end{equation} \end{shrinkeq}
where $x^{(N)}(t)=\frac{1}{N}\sum_{i=1}^Nx_i(t)$ is the average state of the minor agents.
\begin{remark}We remark that the control process and state-average enter both the drift and diffusion terms. This makes our paper different to standard mean-field game (e.g., \cite{WH2017}) or social optimization (e.g., \cite{HCM2012}) literature in which only drift terms are control-dependent.\end{remark}
Let $u(\cdot):=(u_0(\cdot),u_1(\cdot),\cdots,u_N(\cdot))$ be the set of strategies of all $N+1$ agents, $u_{-0}(\cdot):=(u_1(\cdot),\cdots,u_N(\cdot))$ and  $u_{-i}(\cdot):=(u_0(\cdot),\cdots,u_{i-1}(\cdot),u_{i+1}(\cdot),\cdots,u_N(\cdot))$, $0\leq i\leq N$. The centralized admissible strategy set is given by
\begin{shrinkeq}{-1ex}\begin{equation*}
  \begin{aligned}
     &\rr \mathcal U_c\!:=\!\big\{\!u(\cdot)|u(t)\text{ is } \mathcal{F}_t \text{ measurable}, \E \scaleobj{.8}{\int_0^T}\|u(t)\|^2dt\!<\!\infty\!\big\}. \\
  \end{aligned}
\end{equation*}\end{shrinkeq}
{\r Correspondingly}, the feedback decentralized admissible strategy set for the major agent is given by
\begin{shrinkeq}{-1ex}\begin{equation*}\rr
  \resizebox{\linewidth}{!}{$
      \begin{aligned}
         &\mathcal U_0:=\{u_0(\cdot)|u_0(t)\text{ is } \mathcal{F}_t^{W_0}\text{ measurable},  \E \scaleobj{.8}{\int_0^T}\|u_0(t)\|^2dt<\infty\},                                                                         \\
      \end{aligned}
    $}
\end{equation*}\end{shrinkeq}
and the feedback decentralized admissible strategy set for the $i^\text{th}$ minor agent is given by
\begin{shrinkeq}{-1ex}\begin{equation*}\rr
      \begin{aligned}
         & \mathcal U_i:=\{u_i(\cdot)|u_i(t)\text{ is } \mathcal{F}_t^i\text{ measurable},   \E \scaleobj{.8}{\int_0^T}\|u_i(t)\|^2dt<\infty\}.                                                                        \\
      \end{aligned}
\end{equation*}\end{shrinkeq}
For simplicity, define\s
\begin{shrinkeq}{-1ex}\begin{equation*}
  \begin{aligned}
    \mathcal U_{-0}:=\{(u_1(\cdot),\cdots,u_N(\cdot))|u_i(\cdot)\in \mathcal U_i,\ i=1,\cdots,N\}. \\
  \end{aligned}
\end{equation*}\end{shrinkeq}\n
The cost functional for $\mathcal A_0$ is given by
\begin{shrinkeq}{-1ex}\begin{equation}\label{cost-major}
  \resizebox{8.09cm}{!}{$
      \begin{aligned}
         & \hspace{4mm}\mathcal J_0(u_0(\cdot),u_{-0}(\cdot))                                                                \\
         & =\frac{1}{2}\E \!\scaleobj{.8}{\int_0^T}\![\langle Q_0(t)(x_0(t)\- H_0(t)x^{(N)}(t)),x_0(t)\- H_0(t)x^{(N)}(t)\rangle \\
         & \hspace{1.8cm}+\langle R_0(t)u_0(t),u_0(t)\rangle]dt,                                                 \\
      \end{aligned}$}
\end{equation}\end{shrinkeq}
and the cost functional for $\mathcal A_i$, $1\leq i\leq N$, is given by
\begin{shrinkeq}{-1ex}
\begin{equation}\label{cost-minor}\resizebox{8.09cm}{!}{$ \begin{aligned}
         & \hspace{4mm}\mathcal J_i(u_i(\cdot),u_{-i}(\cdot))                                                                                                      \\
         & =\frac{1}{2}\E \scaleobj{.8}{\int_0^T}\Big[\Big\langle Q(t)(x_i(t) \- H(t)x_0(t) \- \hat H(t) x^{(N)}(t)),                                                            \\
         & \hspace{2cm}x_i(t) \- H(t)x_0(t) \- \hat H(t)x^{(N)}(t)\Big\rangle \+ \Big\langle R(t)u_i(t),u_i(t)\Big\rangle\Big]dt. \\
      \end{aligned}$}\end{equation}\end{shrinkeq}
{\dg\begin{remark}
  It is worth pointing out that it brings no essential difficulty to introduce a terminal cost term in \eqref{cost-major} and \eqref{cost-minor}. This will only change the terminal value of the associated Riccati equations.  Thus, for simplicity,  we only consider Lagrange type cost functional here.
\end{remark}}
The aggregate team cost of $N$ minor agents is
\begin{shrinkeq}{-1.5ex}\begin{equation}\label{cost-social optimal}
  \begin{aligned}
    \mathcal J_{soc}^{(N)}(u(\cdot))=\scaleobj{.8}{\sum_{i=1}^N}\mathcal J_i(u_i(\cdot),u_{-i}(\cdot)). \\
  \end{aligned}
\end{equation}\end{shrinkeq}
We impose the following general assumptions, {\r which  are commonly used in LQG models}, on the coefficients:
\begin{itemize}\s
  \item[(H1)] $A_0(\cdot),F_0(\cdot),C_0(\cdot),\widetilde F_0(\cdot),A(\cdot),F(\cdot),C(\cdot),\widetilde F(\cdot),$\\
        $\widetilde G(\cdot) \in L^\infty(0,T;\mathbb R^{n\times n})$,\\ $B_0(\cdot),D_0(\cdot),B(\cdot),D(\cdot)\in L^\infty(0,T;\mathbb R^{n\times m})$.
  \item[(H2)] $Q_0(\cdot)$, $H_0(\cdot)$, $Q(\cdot)$, $H(\cdot)$, $\hat H(\cdot)$ $\in$ $ L^\infty(0,T;\mathbb S^n)$, $R_0(\cdot)$, $R(\cdot)\in L^\infty(0,T;\mathbb S^m)$.\n
\end{itemize}
\begin{remark}
  Under (H1), the system \eqref{state equation-major} and \eqref{state equation-minor} admits a unique strong solution
  $(x_0, \cdots, x_N)$ $\in$ $ L^2_{\F_t}(\Omega;C(0,T;\mathbb{R}^n))$ $\times$ $\cdots$ $\times$ $L^2_{\F_t}(\Omega;C(0,T;\mathbb{R}^n))$ for any given admissible control $(u_0,\cdots,u_N)$ $\in$ $\mathcal U_c$. Under (H2), the cost functionals \eqref{cost-major} and \eqref{cost-minor} are well defined.
\end{remark}
Note that while the coefficients are dependent on the time variable $t$, in what follows, the variable $t$ will usually be suppressed if no confusion would occur.
We propose the following social optimization problem:
\begin{problem}\label{prob1}
Find a decentralized strategy set $\bar u(\cdot)$ = $(\bar u_0(\cdot)$, $\bar u_1(\cdot)$, $\cdots$, $\bar u_N(\cdot))$ where $\bar u_i(\cdot)\in \mathcal U_i$, $0\leq i\leq N$, such that\s
\begin{shrinkeq}{-1ex}\begin{equation}\label{optimal problem}\left\{\begin{aligned}
     & \mathcal J_0(\bar u_0(\cdot),\bar u_{-0}(\cdot))=\inf_{u_0\in\mathcal U_0}\mathcal J_0(u_0(\cdot),\bar u_{-0}(\cdot)),  \\
     & \mathcal J_{soc}^{(N)}(\bar u_0(\cdot),\bar u_{-0}(\cdot))=\inf_{u_{-0}\in\mathcal U_{-0}}\mathcal J_{soc}^{(N)}(\bar u_0(\cdot),u_{-0}(\cdot)).\\
  \end{aligned}\right.\end{equation}\end{shrinkeq}\n
\end{problem}
\section{Auxiliary optimal control problem of the major agent}\label{auxiliary problem-major}
Replacing $x^{(N)}(\cdot)$ of \eqref{state equation-major} and \eqref{cost-major} by $\hat{x}(\cdot)$ which will be determined in Section V, the limiting major agent's state is given by
\begin{shrinkeq}{-1ex}\begin{equation}\label{limit state equation-major}
  \resizebox{8.09cm}{!}{$
  \left\{\begin{aligned}
     & dz_0 \= \big(A_0z_0 \+ B_0v_0 \+ F_0\hat{x}\big)dt \+ \big(C_0z_0 \+ D_0v_0 \+ \widetilde F_0\hat{x}\big)dW_0, \\
     & z_0(0)=\xi_0,                                                                                                                    \\
  \end{aligned}\right.$}
\end{equation}\end{shrinkeq}
and correspondingly the limiting cost functional is
\begin{shrinkeq}{0ex}\begin{equation}\label{limit cost-major}
  \resizebox{8.09cm}{!}{$
      \begin{aligned}
         & J_0(v_0(\cdot))\=\frac{1}{2}\E \!\scaleobj{.8}{\int_0^T}\![\langle Q_0(z_0 \- H_0\hat{x}),z_0 \- H_0\hat{x}\rangle \+ \langle R_0v_0,v_0\rangle]dt. \\
      \end{aligned}
    $}
\end{equation}\end{shrinkeq}
We define the following auxiliary stochastic optimal control problem for major agent:
\begin{problem}\label{prob2}
For agent $\mathcal A_0$, minimize $J_0(u_0(\cdot))$ over {\b$L_{\F^{W_0}_t}^2(0,T;\mathbb{R}^n)$}.
\end{problem}
This is a standard LQ stochastic control problem. For
its solvability, one can introduce the following standard assumption
\begin{description}
  \item[(SA)]  $Q_0(\cdot)\geq0$, $Q(\cdot)\geq0$, $R_0(\cdot)\gg0,$ $R(\cdot)\gg0$.
\end{description}
By \cite[Theorem 4.3]{SLY16}, we have the following result:
\begin{proposition}\label{prop1}
  Under (H1)-(H2) and (SA),  the following Riccati equation
  \begin{shrinkeq}{-1ex}\begin{equation}\label{Riccati-major}
    \resizebox{8.09cm}{!}{$\left\{\begin{aligned}
           & -(P_0B_0+C_0^\top\! P_0D_0)(R_0+D_0^\top\! P_0D_0)^{-1}(B_0^\top\! P_0+D_0^\top\! P_0C_0) \\
           &\+ \dot{P_0}+P_0A_0+A_0^\top\! P_0+C_0^\top\! P_0C_0+Q_0=0, \quad P_0(T)=0,\\
        \end{aligned}\right.$}\end{equation}\end{shrinkeq}
  is strongly regularly solvable, and Problem \ref{prob2} admits a feedback optimal control $\bar v_0\=\Theta_1\bar z_0\+\Theta_2$ where \f
  \begin{shrinkeq}{-1.2ex}\begin{equation}\label{Theta12}
  \left\{\begin{aligned}
       & \Theta_1\=-(R_0\+D_0^\top\! P_0D_0)^{-1}(B_0^\top\! P_0\+D_0^\top\! P_0C_0),                                     \\
       & \Theta_2\=-(R_0\+D_0^\top\! P_0D_0)^{-1}(B_0^\top\phi\+D_0^\top\zeta\+D_0^\top\! P_0\widetilde F_0\hat{x}),
    \end{aligned}\right.
    \end{equation}\end{shrinkeq}\n
  and $\phi$ satisfies
  \begin{shrinkeq}{-1ex}\begin{equation}\label{BSDE-major}
  \resizebox{7.9cm}{!}{$
  \left\{\begin{aligned}
           & d\phi=-\{[A_0^\top-(P_0B_0+C_0^\top\! P_0D_0)(R_0+D_0^\top\! P_0D_0)^{-1}B_0^\top]\phi                                         \\
           & \hspace{7mm}+
          [C_0^\top-(P_0B_0+C_0^\top\! P_0D_0)(R_0+D_0^\top\! P_0D_0)^{-1}D_0^\top]\zeta                                                        \\
           & \hspace{7mm}+\![C_0\!\!^\top\!\-(P_0B_0\+C_0^\top\! P_0D_0)(R_0\+D_0^\top\! P_0D_0)^{-1}D_0^\top]P_0\widetilde F_0\hat{x} \\
           & \hspace{7mm}+P_0F_0\hat{x}-Q_0H_0\hat{x}\}dt+\zeta dW_0, \quad \phi(T)=0.\\
        \end{aligned}\right.$}
        \end{equation}\end{shrinkeq}
  The corresponding optimal state is\s
  \begin{shrinkeq}{-1ex}\begin{equation}\label{Hamilton system-major}\left\{\begin{aligned}
       & d\bar z_0=[(A_0+B_0\Theta_1)\bar z_0+B_0\Theta_2+F_0\hat{x}]dt \\
       &\quad+[(C_0 \+ D_0\Theta_1)
        \bar z_0 \+ D_0\Theta_2 \+ \widetilde F_0\hat{x}]dW_0,\ \rr \z_0(0)=\xi_0.                                                          \\
    \end{aligned}\right.\end{equation}\end{shrinkeq}\n
\end{proposition}
\section{Stochastic optimal control problem for minor agents}
\subsection{Person-by-person optimality}\label{p-b-p optimality}
Let $(\bar u_1(\cdot),\cdots,\bar u_n(\cdot))$ be
centralized optimal strategies of the minor agents. We now perturb
$u_i(\cdot)$ and keep $\bar u_{-i}(\cdot)$=($\bar u_0(\cdot)$, $\bar
u_1(\cdot)$, $\cdots$, $\bar u_{i-1}(\cdot)$, $\bar u_{i+1}(\cdot)$,
$\cdots$, $\bar u_N(\cdot)$) fixed. For $j=1,\cdots,N$, $j\neq i$, denote the
perturbation $\delta u_i(\cdot)=u_i(\cdot)-\bar u_i(\cdot)$, $\delta x_i(\cdot)=x_i(\cdot)-\bar x_i(\cdot)$, $\delta
x_j(\cdot)=x_j(\cdot)-\bar x_j(\cdot)$, $\delta
x^{(N)}=\frac{1}{N}\sum_{j=1}^{N}\delta x_j(\cdot)$,  and
$\delta\mathcal J_j$ is the first variation (Fr\'{e}chet differential)
of $\mathcal J_j$ w.r.t.  $\delta u_j$. Therefore, $\delta x_i$, $\delta x_j$, $\delta x_0$ and $\delta x_{-(0,i)}:=\sum_{j=1,j\neq i}^N\delta x_j$ are given by \s
\begin{shrinkeq}{-1ex}\begin{equation*}
\left\{\begin{aligned}
     & d\delta x_i=(A\delta x_i+B\delta u_i+F\delta x^{(N)})dt \\
     & \hspace{1cm}+(C\delta x_i+D\delta u_i+\widetilde F\delta x^{(N)}+\widetilde G\delta x_0)dW_i,\\
     & d\delta x_j=(A\delta x_j\+F\delta x^{(N)})dt \+(C\delta x_j\+\widetilde F\delta x^{(N)}\+\widetilde G\delta x_0)dW_j, \\
     & d\delta x_0=(A_0\delta x_0+F_0\delta x^{(N)})dt+(C_0\delta x_0+\widetilde F_0\delta x^{(N)})dW_0, \\
     & d\delta x_{-(0,i)}=[A\delta x_{-(0,i)}+F(N-1)\delta x^{(N)}]dt                                    \\
     & \hspace{1.7cm}+ \scaleobj{.8}{\sum_{j\neq i}}(C\delta x_j+\widetilde F\delta x^{(N)}+\widetilde G\delta x_0)dW_j, \\
     & \delta x_i(0)=0,\ \delta x_j(0)=0,\ \delta x_0(0)=0, \  \delta x_{-(0,i)}(0)=0.
  \end{aligned}\right.\end{equation*}\end{shrinkeq} \n
  {\dg By some elementary calculations, we can further obtain
 $\delta\mathcal J_{i}$ of the cost functional of $\mathcal A_i$ as follows\f
\begin{equation*}\label{cost variation-i} \begin{aligned}
     & \delta \mathcal J_{i} \!= \!\E\!\!\int_0^T\!\!\!\! \langle Q(\bar x_i\-\hat H\bar x^{(N)}\!\-H\bar x_0) , \delta x_i\-\hat H\delta x^{(N)}\!\-H\delta x_0\rangle \+ \langle R\bar u_i,\delta u_i\rangle  dt .
  \end{aligned} \end{equation*} \n
For $j\neq i$, $\delta \mathcal J_{j}$ of the cost functional of $\mathcal A_j$ is given by
\begin{equation*}\label{cost variation-j}\resizebox{\linewidth}{!}{$\begin{aligned}
        \delta \mathcal J_{j}=\E\int_0^T\!\!\!\langle Q(\bar x_j-\hat H\bar x^{(N)}-H\bar x_0) , \delta x_j-\hat H\delta x^{(N)}-H\delta x_0\rangle dt.
      \end{aligned}$}\end{equation*}}
We can further obtain
  $\delta \mathcal J_{soc}^{(N)}$, the first variation of the social cost, satisfying
\begin{shrinkeq}{0ex}\begin{equation} \label{cost variation-social cost}
\resizebox{7.9cm}{!}{$
 \begin{aligned}
     & \delta \mathcal J_{soc}^{(N)}\=\E\!\!\int_0^T  \!\! \Big[\langle Q(\bar x_i\-\hat H\bar x^{(N)}\-H\bar x_0),\delta x_i-H\delta x_0-\hat H\delta x^{(N)}\rangle                                         \\
     & +\scaleobj{.8}{\sum_{j\neq i}}\langle Q(\bar x_j\-\hat H\bar x^{(N)}\!\-H\bar x_0),\delta x_j\-H\delta x_0\-\hat H\delta x^{(N)}\rangle \+\langle R\bar u_i,\delta u_i\rangle \Big]dt . \\
  \end{aligned}$} \end{equation}\end{shrinkeq}
Replacing $\bar x^{(N)}$ in \eqref{cost variation-social cost} by $(\bar x^{(N)} - \hat{x}) + \hat{x}$,
\f
  \begin{shrinkeq}{0ex}\begin{equation*} \begin{aligned}
      & \delta \mathcal J_{soc}^{(N)} \= \E\scaleobj{.8}{\int_0^T}  \! [\langle  Q\bar x_i,\delta x_i\rangle\-\langle  Q(\hat H\hat x+H\bar x_0) \+ \hat H Q(\hat x\-\hat H\hat x\-H\bar x_0),\delta x_i\rangle        \\
      &\quad \-\big\langle\hat H Q(\hat x\-\hat H\hat x\-H\bar x_0),\delta x_{\-(0,i)}\big\rangle \-\big\langle HQ(\hat x\-\hat H\hat x\-H\bar x_0),N\delta x_0\big\rangle \\
      &\quad \+ \scaleobj{.8}{\frac{1}{N}\sum_{j\neq i}}\big\langle Q(\bar x_j\-\hat H\hat x\-H\bar x_0),N\delta x_j\big\rangle \+\langle R\bar u_i,\delta u_i\rangle ]dt +  \scaleobj{.8}{\sum_{l=1}^4\varepsilon_{l,i}},
  \end{aligned} \end{equation*}\end{shrinkeq}\n
  where\s
\begin{shrinkeq}{0ex} \begin{equation}\label{e1-4}\left\{\begin{aligned}
     & \varepsilon_{1,i}=\E \scaleobj{.7}{\int_0^T}\langle (Q\hat H-\hat H Q\hat H)(\hat x-\bar x^{(N)}),N\delta x^{(N)}\rangle dt, \\
     & \varepsilon_{2,i}=-\E \scaleobj{.7}{\int_0^T}\langle HQ\hat H(\hat x-\bar x^{(N)}),N\delta x_0\rangle dt,                   \\
     & \varepsilon_{3,i}=\E \scaleobj{.7}{\int_0^T} \langle HQ(\hat x-\bar x^{(N)}),N\delta x_0\rangle dt,                          \\
     & \varepsilon_{4,i}=\E \scaleobj{.7}{\int_0^T} \langle\hat HQ(\hat x-\bar x^{(N)}),N\delta x^{(N)}\rangle dt.
  \end{aligned}\right.\end{equation}\end{shrinkeq}\n
  Introduce the limit processes $(x_0^*, x_j^*, x^{**})$ to replace ($N\delta x_0$, $N\delta x_j$, $\delta x_{-(0,i)}$) by (($N\delta x_0 - x_0^*) + x_0^*$, $(N\delta x_j - x_j^*) + x_j^*$, $(\delta x_{-(0,i)} - x^{**}) + x^{**}$) where
\begin{shrinkeq}{0ex}\begin{equation}\label{limit process of variation}
\resizebox{7.9cm}{!}{$
\left\{\begin{aligned}
     & dx_0^*\=(A_0x_0^*\+F_0\delta x_i\+F_0x^{**})dt \+(C_0x_0^*\+\widetilde F_0\delta x_i\+\widetilde F_0 x^{**})dW_0,                               \\
     & dx_j^*\=(Ax_j^*\+F\delta x_i\+Fx^{**})dt \+\!(Cx_j^* \+ \widetilde F\delta x_i \+ \widetilde F x^{**} \+ \widetilde Gx_0^*)dW_j, \\
     & dx^{**}\=(Ax^{**}\+F\delta x_i\+Fx^{**})dt, \ \  x_0^*(0) \= x_j^*(0) \= x^{**}(0)=0.                                                                                                                           \\
  \end{aligned}\right.$}
  \end{equation}\end{shrinkeq}
Therefore,
\f
\begin{shrinkeq}{0ex}\begin{equation*} \begin{aligned}
    \delta \mathcal J_{soc}^{(N)} = & \E \scaleobj{.8}{\int_0^T}  \Big[\langle Q\bar x_i,\delta x_i\rangle\-\langle Q(\hat H\hat x\+H\bar x_0)  \+\hat H Q(\hat x\-\hat H\hat x\-H\bar x_0),\delta x_i\rangle   \\
                                    & \-\langle \hat HQ(\hat x\-\hat H\hat x\-H\bar x_0),x^{**}\rangle \-\langle HQ(\hat x\-\hat H\hat x\-H\bar x_0), x_0^*\rangle                                   \\
                                    & \+ \scaleobj{.8}{\frac{1}{N}\sum_{j\neq i} } \langle Q(\bar x_j\-\hat H\hat x\-H\bar x_0), x_j^*\rangle \+\langle R\bar u_i,\delta u_i\rangle\Big] dt \+  \scaleobj{.8}{\sum_{l=1}^7 \varepsilon_{1,i}},\\
  \end{aligned}\end{equation*}\end{shrinkeq}\n where \s
\begin{shrinkeq}{0ex}\begin{equation}\label{e5-7}  \left\{\begin{aligned}
     & \varepsilon_{5,i}\=\E \scaleobj{.8}{\int_0^T}\langle\hat H Q(\hat x\-\hat H\hat x\-H\bar x_0),x^{**}\-\delta x_{\-(0,i)}\rangle dt,            \\
     & \varepsilon_{6,i}\=\E \scaleobj{.8}{\int_0^T}\langle HQ(\hat x\-\hat H\hat x\-H\bar x_0),x_0^*\- N\delta x_0\rangle dt,                           \\
     & \varepsilon_{7,i}\=\E \scaleobj{.8}{\int_0^T\frac{1}{N}\sum_{j\neq i}}\langle Q(\bar x_j\-\hat H\hat x\-H\bar x_0),N\delta x_j\- x^*_j\rangle dt.\\
  \end{aligned}\right. \end{equation}\end{shrinkeq}\n
Replacing $\bar x_0$ by $(\bar x_0 - z_0) + z_0$, we have
\s\begin{shrinkeq}{0ex}\begin{equation}\label{variation-conti-1} \begin{aligned}
    \delta \mathcal J_{soc}^{(N)}
    \=  \E \scaleobj{.8}{\int_0^T} & [\langle Q\bar x_i,\delta x_i\rangle\-\langle Q(\hat H\hat x\+Hz_0)\+\hat H Q(\hat x\-\hat H\hat x               \\
                          & -\!Hz_0),\delta x_i\rangle\-\langle \hat HQ(\hat x\-\hat H\hat x\-Hz_0),x^{**}\rangle                                   \\
                          & -\!\langle HQ(\hat x \- \hat H\hat x \- Hz_0), x_0^*\rangle \+ \scaleobj{.8}{\frac{1}{N}\sum_{j\neq i} }\langle Q(\bar x_j              \\
                          & -\!\hat H\hat x \- Hz_0), x_j^*\rangle \+ \langle R\bar u_i,\delta u_i\rangle] dt  \+  \scaleobj{.8}{ \sum_{l=1}^{10}\varepsilon_{l,i}},\\
  \end{aligned} \end{equation}\end{shrinkeq}\n
where\s
\begin{shrinkeq}{-1ex}\begin{equation}\label{e8-10}
 \left\{\begin{aligned}
     & \varepsilon_{8,i}\=\E \!\scaleobj{.8}{\int_0^T}\!\!\langle \hat HQH(\bar x_0\-z_0),x^{**}\rangle dt, \\[-1mm]
     &\hspace{1cm} \+ \E \!\scaleobj{.8}{\int_0^T}\!\!\langle QH(z_0\-\bar x_0)\+\hat HQH(z_0\-\bar x_0),\delta x_i\rangle dt, \\[-1mm]
     & \varepsilon_{9,i}=\E \scaleobj{.8}{\int_0^T}\langle HQH(\bar x_0\-z_0),x_0^*\rangle dt,                                                                        \\[-1mm]
     & \varepsilon_{10,i}=\E \scaleobj{.8}{\int_0^T\frac{1}{N}\sum_{j\neq i}}\langle QH(z_0\-\bar x_0),x_j^*\rangle dt.                                               \\
  \end{aligned}\right.
  \end{equation}\end{shrinkeq}\n
Now we introduce the adjoint equations  $y_1^0$, $y_1^j$ and $y_2$
\s
\begin{shrinkeq}{0ex}\begin{equation*} \left\{\begin{aligned}
     & dy_1^0\=[HQ(\hat x-\hat H\hat x-Hz_0)-A_0^\top\! y_1^0-C_0^\top\beta_1^0                                                                                                       \\
     & \hspace{1cm}-\widetilde G^\top\Eo [\beta_1^{jj}]dt+\beta_1^0dW_0,\quad y_1^0(T)=0,                                                                    \\
     & dy_1^j\=[-Q(\bar x_j-\hat H\hat x-Hz_0)-A^\top\! y_1^j-C^\top\beta_1^{jj}] dt                                                                                              \\
     & \hspace{1cm}\+\beta_1^{jj}dW_j \+ \scaleobj{.8}{\sum_{k\neq j}}\beta_1^{jk}dW_k,\ \ y_1^j(T)=0,\ \ j\=1,\cdots,N,                                                                                                  \\
     & dy_2\=[\hat HQ(\hat x \- \hat H\hat x \- Hz_0) \- F^\top\!\Eo [y_1^j]\-\widetilde F^\top\Eo [\beta_1^{jj}] \\
     & \hspace{1cm}-\!\widetilde F^\top\Eo [\beta_1^{jj}] \- (A \+ F)^\top\! y_2 \- F_0^\top\! y_1^0-\widetilde F_0^\top\beta_1^0]dt                \\
     & \hspace{1cm}+\beta_2^0dW_0,\quad y_2(T)=0,                                                                                                                   \\
  \end{aligned}\right.  \end{equation*}\end{shrinkeq}\n
to replace the terms $x_0^*$, $x_j^*$ and $x^{**}$ respectively,  where \s $\beta_1^0$, $(\beta_1^{j1},\cdots,\beta_1^{jN})$, $\beta_2^0$ \n are the adjoint-states of $y_1^0$, $y_1^j$ and $y_2$ respectively.
 \begin{remark}
    It is important to construct an auxiliary LQG control problem for investigating decentralized control in social optima problem.   Since \eqref{variation-conti-1}, the direct variation decomposition of \s $\delta \mathcal J_{soc}^{(N)}$ \n,  contains $x^{**}$, $x_j^{*}$ and $x_0^{*}$, which are some intermediate variation terms related to the basic variation term $\delta x_i$ indirectly, we have to use a duality procedure (see \cite{HWY2019}) to break away $\delta \mathcal J_{soc}^{(N)}$ from the dependence on $x^{**}$, $x_j^{*}$ and $x_0^{*}$. To this end, we introduce such three adjoint equations.
  \end{remark}
Applying It\^{o}'s formula to $\langle y_1^j,x_j^*\rangle$, $\langle y_2,x^{**}\rangle$ and $\langle y_1^0,x_0^{*}\rangle$, \s
{\dg we have
\small
\begin{align}
    0= & \E\langle y_1^j(T),x_j^*(T)\rangle-\E\langle y_1^j(0),x_j^*(0)\rangle                                                                          \\
  \notag=          & \E\int_0^T\Big[\langle-Q(\bar x_j-\hat H\hat x-Hz_0),x_j^*\rangle+\langle F^\top y_1^j
  +\widetilde F\beta_1^{jj},x^{**}\rangle                                                                                                                                           \\
  \notag           & \hspace{1.2cm}+\langle \widetilde G^\top\beta_1^{jj},x_0^*\rangle + \langle F^\top y_1^j+\widetilde F^\top\beta_1^{jj},\delta x_i\rangle\Big]dt.\\
   =  & \E\langle y_2(T),x^{**}(T)\rangle-\E\langle y_2(0),x^{**}(0)\rangle                                                                            \\
  \notag=          & \E\int_0^T\Big[\langle\hat HQ(\hat x\-\hat H\hat x \- Hz_0) \- F^\top\Eo[y_1^j]                                    \\
  \notag           & \hspace{9mm}\-\widetilde F^\top\Eo[\beta_1^{jj}]\-F_0^\top \!y_1^0\-\widetilde F_0^\top\!\beta_1^0,x^{**}\rangle \+ \langle F^\top\! y_2,\delta x_i\rangle\Big]dt,                                                                                                     \\
   0=  & \E\langle y_1^0(T),x_0^{*}(T)\rangle - \E\langle y_1^0(0),x_0^{*}(0)\rangle                                                                    \\
  \notag=          & \E\int_0^T\Big[\langle HQ(\hat x-\hat H\hat x-Hz_0)-\widetilde G^\top\Eo[\beta_1^{jj}],x_0^*\rangle                              \\
  \notag           & \hspace{1.1cm}+\langle F_0^\top y_1^0+\widetilde F_0^\top\beta_1^0,x^{**}\rangle+\langle F_0^\top y_1^0+\widetilde F_0^\top\beta_1^0,\delta x_i\rangle\Big]dt.
\end{align}\n
Adding to \eqref{variation-conti-1}, we have }
\begin{shrinkeq}{-1ex}
\begin{equation*}\label{variation-conti-2}
\begin{aligned}
        \delta \mathcal J_{soc}^{(N)}
        \=\E\scaleobj{.8}{\int_0^T} & [\langle Q\bar x_i,\delta x_i\rangle \- \langle Q(\hat H\hat x \+ Hz_0)  \+\hat H Q(\hat x     \- \hat H\hat x \- Hz_0)                                          \\
                           &   \- F^\top\! y_2   \- F^\top\E[y_1^j|\mathcal F_t^{W_0}] \- \widetilde F^\top\E[\beta_1^{jj}|\mathcal F_t^{W_0}] \- F_0^\top\! y_1^0 \\
                           & \-\widetilde F_0^\top\beta_1^0,\delta x_i\rangle
        \+\langle R\bar u_i,\delta u_i\rangle ] dt \+  \scaleobj{.8}{\sum_{l=1}^{13}}\varepsilon_l,                                                               \\
      \end{aligned} \end{equation*}\end{shrinkeq} \n
where\f
\begin{shrinkeq}{-1ex}\begin{equation}\label{e11-13} \left\{\begin{aligned}
     & \varepsilon_{11,i}\=\E \scaleobj{.8}{\int_0^T}\langle F^\top( \scaleobj{.8}{\frac{1}{N}\sum_{j\neq i}} y_1^j\-\Eo [y_1^j])                          \\
     & \hspace{1.8cm}+\widetilde F^\top( \scaleobj{.8}{\frac{1}{N}\sum_{j\neq i}} \beta_1^{jj}\-\Eo [\beta_1^{jj}]) ,x^{**}\rangle dt,              \\
     & \varepsilon_{12,i}\=\E \scaleobj{.8}{\int_0^T}\langle\widetilde G^\top( \scaleobj{.8}{\frac{1}{N}\sum_{j\neq i}} \beta_1^{jj}\-\Eo [\beta_1^{jj}]),x_0^*\rangle dt, \\
     & \varepsilon_{13,i}\=\E \scaleobj{.8}{\int_0^T}\langle F^\top( \scaleobj{.8}{\frac{1}{N}\sum_{j\neq i}} y_1^j\-\Eo [y_1^j])                          \\
     & \hspace{1.8cm}+\widetilde F^\top( \scaleobj{.8}{\frac{1}{N}\sum_{j\neq i}} \beta_1^{jj}\-\Eo [\beta_1^{jj}]) ,\delta x_i\rangle dt.
  \end{aligned}\right.\end{equation}\end{shrinkeq}\n
Therefore, considering the case when $N\longrightarrow\infty$, we introduce the first variation of the decentralized auxiliary cost functional  $\delta J_i$ as follows
\begin{shrinkeq}{-1ex}\begin{equation}\label{variation-conti-3}
\resizebox{7.9cm}{!}{$
\begin{aligned}
    &\delta J_i
    =\E\int_0^T  \!\Big[\langle Q\bar x_i,\delta x_i\rangle \- \langle Q(\hat H\hat x \+ Hz_0) \+ \hat H Q(\hat x \- \hat H\hat x \- Hz_0) \\
                       & \-F^\top\! y_2\-F^\top\hat y_1\-\widetilde F^\top\hat\beta_1\-F_0^\top\! y_1^0\-\widetilde F_0^\top\beta_1^0,\delta x_i\rangle  \+\langle R\bar u_i,\delta u_i\rangle\Big] dt.                                                                                            \\
  \end{aligned} $}
  \end{equation}\end{shrinkeq}\n
\begin{remark}
    In \eqref{variation-conti-3}, we ignore $\varepsilon_1,\cdots,\varepsilon_{13}$ and introduce the  first variation of the auxiliary cost functional  $\delta J_i$. Actually, $\varepsilon_1,\cdots,\varepsilon_{13}$ have some order as $\bar x^{(N)} - \hat{x}$, and it is sufficient to conjecture $\|\bar x^{(N)} - \hat{x}\|^2_{L^2}\longrightarrow 0$ when $N\longrightarrow\infty$  by considering the weakly coupled structure of our problem. The rigorous proof will be given in Section VII.
  \end{remark}
\subsection{Decentralized strategy}\label{decentralized strategy}
Motivated by \eqref{variation-conti-3}, one can introduce the following auxiliary problem:
\begin{problem}\label{prob3}
Minimize $J_i(u_i)$ over $u_i\in\mathcal U_i$ where
\begin{shrinkeq}{-1ex}\begin{equation} \label{auxiliary problem-minor}
\resizebox{7.9cm}{!}{$
\left\{\begin{aligned}
     & dz_i=\big(Az_i+Bv_i+F\hat x\big)dt                                                                                              \\
     & \hspace{0.8cm}+\big(Cz_i+Dv_i+\widetilde  F\hat x+\widetilde Gz_0\big)dW_i,\hspace{2mm} z_i(0)=z,                                 \\
     & J_i(v_i) \= \frac{1}{2}\E \scaleobj{.75}{\int_0^T}[\langle Q z_i, z_i\rangle\-2\langle S, z_i\rangle\+\langle R v_i, v_i\rangle] dt,\\
     &S=  Q(\hat H\hat x \+ Hz_0) + \hat H Q(\hat x \- \hat H\hat x \- Hz_0) \- F^\top\! y_2  \\
     &\qquad\- F^\top\hat y_1 -\!\widetilde F^\top\hat\beta_1 \- F_0^\top\! y_1^0 \- \widetilde F_0^\top\beta_1^0.
  \end{aligned}\right.$}
  \end{equation}\end{shrinkeq}\n
\end{problem}
The mean-field terms $\hat x$, $z_0$, $y_2$, $\hat y_1$, $\hat\beta_1$, $y_1^0$, $\beta_1^0$ will be determined by the CC system in Section V.
From \cite{SLY16}, we have the following result:
\begin{proposition}\label{prop2}
  Under (H1)-(H2) and (SA), the following Riccati equation\s
  \begin{shrinkeq}{-1ex}\begin{equation}\label{Riccati-minor} \left\{\begin{aligned}
       & \dot{P}\+PA\+A^\top\! P\+C^\top\! PC\+Q   \-(PB \+ C^\top\! PD)                       \\
       &\times\! (R \+ D^\top\! PD)^{-1}(B^\top\! P \+ D^\top\! PC)\=0,\quad P(T)\=0,
    \end{aligned}\right. \end{equation}\end{shrinkeq}\n
  is  strongly regularly solvable, and Problem \ref{prob3} admits a  feedback optimal control $\bar v_i=\Lambda_1\bar x_i+\Lambda_2$ where \s
  \begin{shrinkeq}{-.5ex}\begin{equation}\label{Lambda1,2} \left\{\begin{aligned}
       & \Lambda_1=-(R+D^\top\! PD)^{-1}(B^\top\! P+D^\top\! PC),                                                     \\
       & \Lambda_2\= -(R \+ D^\top\! PD)^{-1}(B^\top\varphi \+ D^\top\eta \+ D^\top\! P(\widetilde  F\hat x \+ \widetilde Gz_0)),
    \end{aligned}\right. \end{equation}\end{shrinkeq}\n
  and $(\varphi,\eta)$ satisfies\s
  \begin{shrinkeq}{-1ex}\begin{equation}\label{BSDE-minor} \left\{\begin{aligned}
       & d\varphi\=-\{[A^\top-(PB+C^\top\! PD)(R+D^\top\! PD)^{-1}B^\top]\varphi                   \\
       & \hspace{7mm}+[C^\top-(PB+C^\top\! PD)(R+D^\top\! PD)^{-1}D^\top]\eta                        \\
       & \hspace{7mm}+[(PB+C^\top\! PD)(R+D^\top\! PD)^{-1}D^\top-C^\top]                            \\
       & \hspace{7mm}\times\! P(\widetilde  F\hat x \+ \widetilde Gz_0) \+ PF\hat x \- S\}dt \+ \eta dW_0, \ \varphi(T) \= 0.
    \end{aligned}\right. \end{equation}\end{shrinkeq}\n
\end{proposition}
  \section{Consistency condition}
Because of the symmetric and decentralized character, we only need a generic Brownian motion (still denoted by $W_1$) which is independent of $W_0$ to characterize the CC system.
\begin{proposition}
  The undetermined quantities in Problem \ref{prob2}, \ref{prob3} can be determined by $(\hat x$, $z_0$, $y_1^0$, $\beta_1^0$, $\hat y_1$, $\hat\beta_1$, $y_2)$=$(\Eo [z]$, $z_0$, $\check y_0$, $\check\beta_0$, $\Eo [\check y_1]$, $\Eo [\check\beta_1^1]$, $\check y_2)$, where $(z$, $z_0$, $\check y_0$, $\check\beta_0$, $\check y_1$, $\check\beta_1^1$, $\check y_2)$ is the solution of the following {\g MF-FBSDEs}: \scriptsize
  \begin{shrinkeq}{-1ex}\begin{equation}\label{CC}\left\{\begin{aligned}
       & d z_0=[(A_0\-B_0\mathcal R_0^{\-1}\cP_0)z_0\-B_0\mathcal R_0^{\-1}B_0^\top\check\phi  \- B_0\mathcal R_0^{\-1}D_0^\top\check\zeta                                                                                                        \\
       & \hspace{10mm}+(F_0\-B_0\mathcal R_0^{\-1}D_0^\top P_0\widetilde F_0)\Eo[z]]dt\+[(C_0\-D_0\mathcal R_0^{\-1}\cP_0)z_0                                                                                                                                             \\
       & \hspace{10mm}\-D_0\mathcal R_0^{\-1}B_0^\top\check\phi\-D_0\mathcal R_0^{\-1}D_0^\top\check\zeta                                                                                                                                       \\
       & \hspace{10mm}+(\widetilde F_0
      \-D_0\mathcal R_0^{\-1}D_0^\top P_0\widetilde F_0)
      \Eo[z]]dW_0,                                                                                                                                                                                           \\
       & dz=[(A\-B\mathcal R^{\-1}\cP)z\-B\mathcal R^{\-1}D^\top P\widetilde Gz_0   \-B\mathcal R^{\-1}\! B^\top\!\check\varphi  \-B\mathcal R^{\-1}\! D^\top\!\!\check\eta                                                                                                                             \\
       & \hspace{10mm}\+\!(\!F\-B\mathcal R^{\-1}\!D^\top\!\! P\widetilde F )\Eo[z] ]dt \+[(C\-D\mathcal R^{\-1}\cP)z               \\
       & \hspace{10mm}\+(\!\widetilde G\-D\mathcal R^{\-1}\!D^\top\!\widetilde G)z_0  \-D\mathcal R^{\-1}B^{\!\top}\!\check\varphi\-D\mathcal R^{\-1}D^\top\!\check\eta                                                                \\
       & \hspace{10mm}+(\widetilde F\-D\mathcal R^{\-1}D^\top P\widetilde  F)\Eo[z]]dW_1,                                                                                                              \\
       & d\check y_0\=[-HQHz_0 \+ HQ(I\-\hat H)\Eo[z]\-A_0^{\!\top}\! \check y_0 \-C_0^{\!\top}\!\check\beta_0\-\widetilde G^\top\Eo[\check\beta_1^1]]dt\\
       &\hspace{10mm}+\check\beta_0dW_0,                                                                                                          \\
       & d\check{y}_1\!=[QHz_0\-Qz\+Q\hat H\Eo[z] \-A^{\!\top}\!\check{y}_1 \-C^{\!\top}\!\check{\beta}_1^1]dt \+\check\beta_1^0dW_0\+\check{\beta}_1^1dW_1,                                 \\
       & d\check{y}_2=[-\hat HQHz_0+\hat HQ(I\-\hat H)
      \Eo[z]\- F^{\!\top}\!\Eo[\check{y}_1]                                                                                                                                                   \\
       & \hspace{10mm}\-\!\widetilde F^\top\!\Eo[\check{\beta}_1^1]\-(\!A\+F)\!^\top \!\check y_2\-F_0\!\!^\top\!\check y_0\-\widetilde F_0\!\!^\top\!\check\beta_0]dt\+\check\beta_2dW_0, \\
       & d\check\phi=\-[(A_0\!\!\!\!^\top\!\-\cPTm\mathcal R_0^{\-1}\!B_0^\top)\check\phi \+ (C_0\!\!\!^\top\!\-\cPTm\mathcal R_0^{\-1}\!D_0^\top)\check\zeta \-(\cPTm\mathcal R_0^{\-1}\!D_0^\top\!\!\-C_0\!\!\!^\top)                                                      \\
       & \hspace{10mm}\times\!P_0\widetilde F_0\Eo[z]     \+P_0F_0\Eo[z] \- Q_0H_0\Eo[z]]dt \+\check\zeta dW_0,                                                                                                  \\
       & d\check{\varphi}=\-[(A\!^\top\!\!\-\cPT\mathcal R^{\-1}\!B\!^\top)\check{\varphi} \+ (C\!^\top\!\!\-\cPT\mathcal R^{\-1}\!D\!^\top)\check{\eta}\+(\cPT\mathcal R^{\-1}\!D\!^\top\!\!\-C^\top\!)P\widetilde Gz_0                                                                                                                   \\
       & \hspace{8mm}\+(PF \+(\cPT\mathcal R^{\-1}\!D\!\!^\top \!\!\- C\!^\top)P\widetilde  F)\Eo[z]\-S ]dt\+\check\eta dW_0,    \\
       & z_0(0)=\xi_0,\quad z(0)=\xi,\quad \check y_0(T)=0,\quad \check{y}_1(T)=0,                                                                                                                                                          \\
       & \check{y}_2(T)=0,\quad \check\phi(T)=0,\quad
      \check{\varphi}(T)=0,
    \end{aligned}\right.\end{equation}\end{shrinkeq}\n
  with\s
   \begin{equation}
   \left\{\begin{aligned}
   &\cP:=\!B^\top\! P\+D^\top\! PC,\quad \cP_0:=\!B_0^\top\! P_0\+D_0^\top\! P_0C_0\\
   &\mathcal R:=R+D^\top\! PD,\quad \mathcal R_0:=R_0+D_0^\top\! P_0D_0.
   \end{aligned}\right.
   \end{equation}\n
\end{proposition}
  Define \f $X=(z_0^\top,z^\top)^\top$,
$Y=(\check{y}_0^\top, \check{y}_1^\top, \check{y}_2^\top, \check\phi^\top, \check{\varphi}^\top)^\top$, $Z_1$=$(\check{\beta}_0^\top$, $\check{\beta}_1^{0\top}$, $\check{\beta}_2^\top$, $\check\zeta^\top$, $\check \eta^\top)^\top$, $Z_2$=$(0^\top$, $\check{\beta}_1^{1\top}$, $0^\top$, $0^\top$, $0^\top)^\top$, $Z=\left(Z_1, Z_2\right)$ \n and  \f $W=({W_0^\top},{W_1^\top})^\top$, \n \eqref{CC} take the following form:\f
  \begin{shrinkeq}{-1ex}\begin{equation}\label{CC-equivalent form} \left\{\begin{aligned}
     & dX=[\mathbb A_1 X+\bar{\mathbb A}_1\E[X|\mathcal F_t^{W_0}]+\mathbb B_1Y+\mathbb F_1Z_1]dt                           \\
     & \hspace{7mm}+[\mathbb C_1^0X+\bar{\mathbb C}_1^0\E[X|\mathcal F_t^{W_0}]+\mathbb D_1^0 Y+\mathbb F_1^0Z_1] dW_0     \\
     & \hspace{7mm}+[\mathbb C_1^1X+\bar{\mathbb C}_1^1\E[X|\mathcal F_t^{W_0}]+\mathbb D_1^1Y+\mathbb F_1^1Z_1]dW_1,      \\
     & dY=[\mathbb A_2X+\bar{\mathbb A}_2\E[X|\mathcal F_t^{W_0}]+\mathbb B_2Y+\bar{\mathbb B}_2\E[Y|\mathcal F_t^{W_0}] \\
     & \hspace{7mm}+\mathbb C_2Z_1 \+ \tilde{\mathbb C}_2Z_2 \+ \bar{\mathbb C}_2\E[ Z_2|\mathcal F_t^{W_0}]
    ]dt \+ Z_1dW_0 \+ Z_2dW_1,                                                                                                               \\
     & X(0)=(\xi_0^\top,\xi^\top)^\top,\quad Y(T)=(0^\top,0^\top,0^\top,0^\top,0^\top)^\top,
  \end{aligned}\right.\end{equation}\end{shrinkeq}
  \n
where \\
  \f $\mathbb A_1=\left(
        \begin{smallmatrix}
            A_0-B_0\mathcal R_0^{\-1}\!\cP_0 & 0 \\
            -B\mathcal R^{\-1}\!D^{\!\top}\! P\widetilde G & A-B\mathcal R^{\-1}\!\cP \\
          \end{smallmatrix}
        \right)$, $\bar{\mathbb A}_1=\left(
        \begin{smallmatrix}
            0 & F_0-B_0\mathcal R_0^{\-1}\!D_0^{\!\top}\! P_0\widetilde F_0 \\
            0 & F-B\mathcal RD^{\!\top}\! P\widetilde  F \\
          \end{smallmatrix}
        \right)$, $\mathbb B_1\=\left(
        \begin{smallmatrix}
            0 & 0 & 0 & -B_0\mathcal R_0^{\-1}\!B_0^{\!\top}\! & 0 \\
            0 & 0 & 0 & 0 & -B\mathcal R B^{\!\top}\! \\
          \end{smallmatrix}
        \right)$,
        $\mathbb F_1\=\left(
        \begin{smallmatrix}
            0 & 0 & 0 & -B_0\mathcal R_0^{\-1}\! D_0^{\!\top}\! &0 \\
            0 & 0 & 0 & 0& -B_0\mathcal R_0^{\-1}\! D^{\!\top}\! \\
          \end{smallmatrix}
        \right)$, $\mathbb F_1^0=\left(
        \begin{smallmatrix}
            0 & 0 & 0 & -D_0\mathcal R_0^{\-1}\! D_0^{\!\top}\! &0 \\
            0 & 0 & 0 & 0& 0  \\
          \end{smallmatrix}
        \right)$,
        $\bar{\mathbb C}_1^0=\left(
        \begin{smallmatrix}
            0 & \widetilde F_0
            -D_0\mathcal R_0^{\-1}\!D_0^{\!\top}\! P_0\widetilde F_0 \\
            0 & 0 \\
          \end{smallmatrix}
        \right)$, $\mathbb C_1^0=\left(
        \begin{smallmatrix}
            C_0-D_0\mathcal R_0^{\-1}\!\cPm & 0 \\
            0 & 0 \\
          \end{smallmatrix}
        \right)$,
        $\mathbb C_1^1=\left(
        \begin{smallmatrix}
            0 & 0 \\
            \widetilde G-D\mathcal R^{\-1}\!D^{\!\top}\!\widetilde G & C-D\mathcal R^{\-1}\!\cP \\
          \end{smallmatrix}
        \right)$, $\bar{\mathbb C}_1^1=\left(
        \begin{smallmatrix}
            0 & 0 \\
            0 & \widetilde F-D\mathcal R^{\-1}\!D^{\!\top}\! P\widetilde  F \\
          \end{smallmatrix}
        \right)$, $\mathbb D_1^0=\left(
        \begin{smallmatrix}
            0 & 0 & 0 & -D_0\mathcal R_0^{\-1}\!B_0^{\!\top}\! & 0 \\
            0 & 0 & 0 & 0 & 0 \\
          \end{smallmatrix}
        \right)$, $\mathbb D_1^1=\left(
        \begin{smallmatrix}
            0 & 0 & 0 & 0 & 0 \\
            0 &0 & 0 & 0 & -D\mathcal R^{\-1}\!B^{\!\top}\! \\
          \end{smallmatrix}
        \right)$, $\mathbb F_1^1 = \left(
        \begin{smallmatrix}
            0 & 0 & 0 & 0 &0 \\
            0 & 0 & 0 & 0& -D\mathcal R^{\-1}\!D^{\!\top}\! \\
          \end{smallmatrix}
        \right)$, $\mathbb A_2=\left(
        \begin{smallmatrix}
            -HQH&0\\
            QH & -Q \\
            -\hat HQH&0 \\
            0 & 0 \\
            -(\cP^{\!\top}\!\mathcal R^{\-1}\!D^{\!\top}\!-C^{\!\top}\!)P\widetilde G + QH  -\hat H Q H&0\\
          \end{smallmatrix}
        \right)$, $\bar{\mathbb A}_2=\left(
        \begin{smallmatrix}
            0& HQ(I-\hat H)\\
            0 & Q\hat H\\
            0 &\hat HQ(I-\hat H) \\
            0&\bar{\mathbb A}_2'\\
            0&\bar{\mathbb A}_2''\\
          \end{smallmatrix}
        \right)$, $\mathbb B_2=\left(
        \begin{smallmatrix}
            -A_0^{\!\top}\!&0&0&0&0\\
            0&  -A^{\!\top}\! & 0 & 0&0\\
            -F_0^{\!\top}\!& 0 & -(A+F)^{\!\top}\! & 0&0\\
            0&0&0&\mathbb B_2'&0\\
            -F_0^{\!\top}\!&0&  -F^{\!\top}\! & 0 & \mathbb B_2''\\
          \end{smallmatrix}
        \right)$,  $\bar{\mathbb A}_2'=[\cPTm\mathcal
        R_0^{\-1}\!D_0^{\!\top}\!-C_0^\top]P_0\widetilde
        F_0 \- P_0F_0 \+ Q_0H_0$, $\bar{\mathbb A}_2'' \= \-F\+[\cPT\mathcal
        R^{\-1}D^{\!\top}\!\-C^\top]P\widetilde  F\+ \hat H Q \+ Q\hat H \-
        \hat H Q\hat H$, $\mathbb B_2' = \-A_0^{\!\top}\!+\cPTm\mathcal
        R_0^{\-1}\!B_0^{\!\top}$,  $\mathbb B_2'' =
        -A^{\!\top}\!+\cPT\mathcal R^{\-1}\!B^{\!\top}\!$, $\bar{\mathbb
        B}_2=\left(
        \begin{smallmatrix}
            0&0&0&0&0\\
            0&0 & 0 & 0&0\\
            0& - F^{\!\top}\! & 0 & 0&0\\
            0& 0 & 0 & 0&0\\
            0& - F^{\!\top}\! & 0 & 0&0\\
          \end{smallmatrix}
        \right)$, $\bar{\mathbb C}_2\=\left(
        \begin{smallmatrix}
            0&  -\widetilde G^\top & 0 & 0 & 0\\
            0& 0 & 0 & 0 & 0\\
            0&-\widetilde F^{\!\top}\! & 0& 0 & 0\\
            0& 0 & 0 & 0& 0 \\
            0& -\widetilde F^{\!\top}\! & 0 & 0& 0 \\
          \end{smallmatrix}
        \right)$,    $\mathbb C_2'\=-(C_0^{\!\top}\!-\cPTm\mathcal R_0^{\-1}\!D_0^{\!\top}\!)$,  \!\!\! $\mathbb C_2'' \= -(C^{\!\top} \- \cPT\mathcal R^{\-1}\!D^{\!\top}\!)$.  \n
  Next we  use discounting method to study the global solvability of FBSDEs \eqref{CC-equivalent form}.
  {\dg To start, we first give some results for general nonlinear forward-backward system:\s
\begin{equation}\label{nonlinear FBSDE} \left\{\begin{aligned}
     & dX(t)\!=\!b\Big(t,X(t),\Eo[X(t)],Y(t),Z(t)\Big)dt                              \\
     & \hspace{1.2cm}+\!\sigma\Big(t,X(t),\Eo[X(t)], \!Y(t),Z(t)\Big)dW(t),           \\
     & dY(t)\!=\!-f\Big(\!t,X(t),\Eo[X(t)],Y(t),\E[Y(t)|\F_t^{W_t}], \\
     & \hspace{2.3cm}Z(t),\Eo[Z(t)]\Big)dt+Z(t)dW(t),                                 \\
     & X(0)=x,\qquad Y(T)=0,
  \end{aligned}\right.\end{equation}\n
where $W = \binom{W_0}{W_1}$, and the coefficients satisfy the following conditions:
\begin{description}
  \item[(A1)] There exist $\rho_1,\rho_2\in\mathbb R$ and positive constants $k_i,i=1,\cdots,12$ such that
        for all $t,x,\bar x,y,\bar y,z,\bar z$, a.s.,
        \s\begin{enumerate}[leftmargin=*]
          \item $\langle b(t,x_1,\bar x,y,z)- b(t,x_2,\bar x,y,z),x_1-x_2\rangle\leq \rho_1\|x_1-x_2\|^2,$
          \item $\|b(t,x,\bar x_1,y_1,z_1)\-b(t,x,\bar x_2,y_2,z_2)\|\!\leq\! k_1\|\bar x_1\-\bar x_2\|+k_2\|y_1-y_2\|+k_3\|z_1-z_2\|,$
          \item $\langle f(t,x,\bar x,y_1,\bar y,z,\bar z) - f(t,x,\bar x,y_2,\bar y,z,\bar z),$\\
                $y_1-y_2\rangle\leq\rho_2\|y_1-y_2\|^2,$
          \item $\|f(t,x_1,\bar x_1,y,\bar y_1,z_1,\bar z_1)-f(t,x_2,\bar x_2,y,\bar y_2,z_2,$\\$\bar z_2)\|\leq k_4\|x_1-x_2\|+k_5\|\bar x_1-\bar x_2\|+k_6\|\bar y_1-\bar y_2\|$\\
                $+k_7\|z_1-z_2\|+k_8\|\bar z_1-\bar z_2\|$,
          \item $\|\sigma(t,x_1,\bar x_1,y_1,z_1)-\sigma(t,x_2,\bar x_2,y_2,z_2)\|^2\leq$\\
                $k_9^2\|x_1-x_2\|^2+k_{10}^2\|\bar x_1-\bar x_2\|^2+k_{11}^2\|y_1-y_2\|^2+k_{12}^2\|z_1-z_2\|^2.$
        \end{enumerate}\n
  \item[(A2)]\s
        \begin{equation}
          \begin{aligned}
            \hspace{-8mm}\E\!\!\int_0^T\!\! & \|b(t,0,0,0,0)\|^2 \!\+ \|\sigma(t,0,0,0,0)\|^2 \!\+ \|f(t,0,0,0,0,0,0)\|^2 dt \!<\!\! \infty.
          \end{aligned}
        \end{equation}\n
\end{description}

Similar to \cite{HHN2018}, we have the following result of solvability of \eqref{CC-equivalent form}. For the readers' convenience, we give the proof in the appendix.
\begin{theorem}\label{discounting}
  Suppose (A1) and (A2) hold. There exists a constant $\delta_1>0$ depending on $\rho_1$, $\rho_2$, $T$, $k_i$, $i$ = $1$, $6$, $7$, $8$, $9$, $10$ such that if $k_i\in[0,\delta_1)$, $i$ = $2$, $3$, $4$, $5$, $11$, $12$, FBSDE \eqref{nonlinear FBSDE} admits a unique adapted solution $(X,Y,Z)\in L^2_{\F}(0,T;\mathbb R^n)\times L^2_{\F}(0,T;\mathbb R^m)\times L^2_{\F}(0,T;\mathbb R^{m\times d}).$ Furthermore, if $2\rho_1+2\rho_2<-2k_1-2k_6-2k_7^2-2k_8^2-k_9^2-k_{10}^2$,
  there exists a constant $\delta_2>0$ depending on $\rho_1,\rho_2,k_i,i=1,6,7,8,9,10$ such that if $k_i\in[0,\delta_2)$, $i=2,3,4,5,11,12$, FBSDE  \eqref{nonlinear FBSDE} admits a unique adapted solution $(X,Y,Z)\in L^2_{\F}(0,T;\mathbb R^n)\times L^2_{\F}(0,T;\mathbb R^m)\times L^2_{\F}(0,T;\mathbb R^{m\times d}).$
\end{theorem}
}
  Let $\rho_1^*$ and $\rho_2^*$ be the largest eigenvalue of $\frac{1}{2}(\mathbb A_1+\mathbb A_1^\top)$ and $\frac{1}{2}(\mathbb B_2+\mathbb B_2^\top)$ respectively. {\dg Comparing \eqref{nonlinear FBSDE} with \eqref{CC-equivalent form}, we can check that the parameters of (A1) can be chosen as follows:} \f
  \begin{shrinkeq}{-1ex}\begin{equation*} \begin{aligned}
     & k_1=\|\bar{\mathbb A}_1\|,\quad k_2=\|\mathbb B_1\|,\quad k_3 = \|F_1\|,\quad k_4=\|\mathbb A_2\|,\quad  k_5=\|\bar{\mathbb A}_2\|,                                             \\
     &  k_6=\|\bar{\mathbb B}_2\|,\quad  k_7=\|\mathbb{C}_2\| + \|\tilde{\mathbb{C}}_2\|,\quad k_8=\|\bar{\mathbb C}_2\|,\quad k_9=\|\mathbb{C}_1^0\| + \|\mathbb{C}_1^1\|, \\
     &  k_{10}=\|\bar{\mathbb C}_1^0\| + \|\bar{\mathbb C}_1^1\|,\quad k_{11}=\|\mathbb{D}_1^0\| + \|\mathbb{D}_1^1\|,\quad k_{12}=\|F_1^0\| + \|F_1^1\|.
  \end{aligned}\end{equation*}\end{shrinkeq}\n
  Now we introduce the following assumption:
\begin{description}
  \small\item[(H3)]
        $2\rho_1^*+2\rho_2^*< - 2\|\bar{\mathbb A}_1\| - 2\|\bar{\mathbb B}_2\| - 2(\|\mathbb{C}_2\| + \|\tilde{\mathbb{C}}_2\|)^2 - 2\|\bar{\mathbb C}_2\|^2 - (\|\mathbb{C}_1^0\| + \|\mathbb{C}_1^1\|)^2 - (\|\bar{\mathbb C}_1^0\| + \|\bar{\mathbb C}_1^1\|)^2.$
\end{description}
We have the following result:
\begin{proposition}\label{prop4}
  Under (H1)-(H3),
  there exists a constant $\delta_3>0$ depending on $\rho_1^*$, $\rho_2^*$, $k_i$, $i=1,6,7,8,9,10,$ such that if $k_i\in[0,\delta_3)$, $i=2,3,4,5,11,12,$ {\g FBSDEs}  \eqref{CC-equivalent form} admits a unique adapted solution $(X,Y,Z)\in L^2_{\F}(0,T;\mathbb R^n)\times L^2_{\F}(0,T;\mathbb R^m)\times L^2_{\F}(0,T;\mathbb R^{m\times d}).$
\end{proposition}
 In what follows, we give an example to show how exactly such conditions can be applied.
  \begin{remark}
    For $\varepsilon > 0$, let  $\rho_1 \= \frac{k_2}{\varepsilon}$, $\rho_2 \= \frac{k_3}{\varepsilon}$, $\rho_3 \= \frac{k_4}{\varepsilon}$, $\rho_4 \= \frac{k_5}{\varepsilon}$, $\rho_5 \= \frac{k_7}{2k_7^2 + 2\varepsilon}$, $\rho_6 \= \frac{k_8}{2k_8^2 + 2\varepsilon}$, $d \= -2k_1 - 2k_6 - 2k_7^2 - 2k_8^2 - k_9^2 - k_10^2 - 2\rho^*_1 - 2\rho^*_2 - 4\varepsilon$, $\bar{\rho}_1 \= \bar{\rho}_2 \= \frac{d}{2}$, $\theta \= \big(\frac{1}{\bar{\rho}_2} + \frac{1}{1 - k_7\rho_5 - k_8\rho_{10}}\big)\big(\frac{1}{\bar{\rho}_1}\big) \= \big(\frac{2}{d} + \frac{(k_7^2 +\varepsilon)(k_8^2+\varepsilon)}{\varepsilon^2}\big) \big(\frac{2}{d}\big)$.
    If $d > 0$, $\theta\big(\frac{k_4^2}{\varepsilon} + \frac{k_5^2}{\varepsilon}\big) < 1$, $\theta\big(\frac{k_2^2}{\varepsilon} + k_{11}^2\big) < 1$, $\theta\big(\frac{k_3^2}{\varepsilon} + k_{12}^2\big) < 1$, then \eqref{CC-equivalent form} admits a unique solution.
  \end{remark}
  Thus, via Propositions \ref{prop1}, \ref{prop2}, \ref{prop4}, we can establish the following procedure to calculate the mean-field strategy.
  \begin{itemize}
    \item Under (H1)-(H3) and (SA), each agent can calculate CC system \eqref{CC} and obtain $(z$, $z_0$, $\check y_0$, $\check\beta_0$, $\check y_1$, $\check\beta_1^1$, $\check y_2)$. Then by taking expectation, the mean-field terms can be obtained by $(\hat x$, $z_0$, $y_1^0$, $\beta_1^0$, $\hat y_1$, $\hat\beta_1$, $y_2)$=$(\Eo [z]$, $z_0$, $\check y_0$, $\check\beta_0$, $\Eo [\check y_1]$, $\Eo [\check\beta_1^1]$, $\check y_2)$.
    \item With $(\hat x$, $z_0$, $y_1^0$, $\beta_1^0$, $\hat y_1$, $\hat\beta_1$, $y_2)$, the agent can solve Riccati equations \eqref{Riccati-major}, \eqref{Riccati-minor} and BSDEs \eqref{BSDE-major}, \eqref{BSDE-minor} to obtain $P_0$, $P$, $\phi$, $\varphi$.
    \item With  $P_0$, $P$, $\phi$, $\varphi$, the agent can obtain $(\Theta_1,\Theta_2)$ and $(\Lambda_1,\Lambda_2)$ by \eqref{Theta12} and \eqref{Lambda1,2} respectively. Then the mean-field  decentralized strategies are given by $ \tilde u_0$ $=$ $\Theta_1\z_0$ $+$ $\Theta_2$, $ \tilde u_i$ $=$ $\Lambda_1\z_i$ $+$ $\Lambda_2$, for $i$ $=$ $1$, $\cdots$, $N$, where  $\z_0$ and $\z_i$  satisfy\s
          \begin{shrinkeq}{0ex}\begin{equation}\label{realized state} \left\{\begin{aligned}
               & d\z_0=(\left(A_0\+B_0\Theta_1\right)\z_0 \+ B\Theta_2 \+F_0\Eo[z])dt \\
               & \hspace{8mm}\+(\left(C_0\+D_0\Theta_1\right)\z_0 \+D_0\Theta_2 \+\widetilde F_0\Eo[z])dW_0, \\
               & d\z_i=(\left(A\+B\Lambda_1\right)\z_i \+ B\Lambda_2 \+F\Eo[z])dt       \\
               & \hspace{8mm}\+(\left(C
                 \+D\Lambda_1\right)\z_i \+ D\Lambda_2\+\widetilde F\Eo[z]\+\widetilde G\z_0)dW_i, \\
               & \z_0(0)=\xi_0,\quad \z_i(0)=\xi, \quad 1\leq i\leq N
            \end{aligned}\right. \end{equation}\end{shrinkeq}\n
  \end{itemize}
  Through the discussion above, the mean-field decentralized strategies are characterized. In what follows, we will show some special cases and illustrate the relation between our research and the existing literature.
\section{Special case}
In this section, we compare our result with standard LQG control problem \cite{SLY16} and LQG social optima \cite{QHX2020}.

When there involves no  minor agent, this problem will reduce to a standard LQG optimal control problem. By letting $A$ = $B$ = $F$ = $C$ =$D$ = $\tilde{F}$ = $\tilde{G}$ = $Q$ = $H$ = $\hat{H}$ = $R$ = 0, we have $x_i^{(N)} = x_i\equiv 0$, and also by \eqref{CC}, we have $z = \hat{x} \equiv 0$. Then by \eqref{Riccati-major} and \eqref{BSDE-major}, we have
\begin{shrinkeq}{-1ex}\begin{equation*} \resizebox{\linewidth}{!}{$\left\{\begin{aligned}
         & -(P_0B_0+C_0^\top\! P_0D_0)(R_0+D_0^\top\! P_0D_0)^{-1}(B_0^\top\! P_0+D_0^\top\! P_0C_0)+ \\
         & \dot{P_0}+P_0A_0+A_0^\top\! P_0+C_0^\top\! P_0C_0+Q_0=0,\quad P_0(T)=0,
      \end{aligned}\right.$}\end{equation*}\end{shrinkeq}
and $\phi\equiv 0$.
Thus, $(\Theta_1, \Theta_2)$ takes the following form
\begin{shrinkeq}{-1ex}\begin{equation*} \begin{aligned}
     & \Theta_1=-(R_0+D_0^\top\! P_0D_0)^{-1}(B_0^\top\! P_0+D_0^\top\! P_0C_0), \quad \Theta_2=0.
  \end{aligned} \end{equation*}\end{shrinkeq}
Such result is consistent with \cite[Theorem 4.3]{SLY16}.

On the other hand, when there involves no major agent, this problem will reduce to a problem of mean-field control in  social optima. By letting $A_0$ = $B_0$ = $F_0$ = $C_0$ =$D_0$ = $\tilde{F}_0$ = $\tilde{G}_0$ = $Q_0$ = $H_0$ = $\hat{H}_0$ = $R_0$ = 0, we have $x_0\equiv 0$.
The CC system becomes:\f
\begin{shrinkeq}{-1ex}\begin{equation*}  \left\{\begin{aligned}
     & dz\=[(A\-B\mathcal R^{\-1}\cP)z \- B\mathcal R^{\-1} B^\top\check\varphi \+(F\-B\mathcal R^{\-1}D^\top P\widetilde F ) \E[z] ]dt                                                                               \\
     & \hspace{7mm}\+[(C\-D\mathcal R^{\-1}\cP)z \- D\mathcal R^{\-1}B^{\!\top}\!\!\check\varphi \+(\widetilde F\-D\mathcal R^{\-1}D^{\!\top}\!\! P\widetilde  F)\E[z]]dW_1(t),                                                                                 \\
     & d\check y_0\=[HQ(I\-\hat H)\E[z]\-\widetilde G^\top\E[\check\beta_1^1]]dt,                                                                                           \\
     & d\check{y}_1\=[-Qz\+Q\hat H\E[z]\-A^\top\check{y}_1\-C^\top\check{\beta}_1^1]dt \+\check{\beta}_1^1dW_1(t),                                                              \\
     & d\check{y}_2\=[ \hat HQ(I\-\hat H)
      \E[z]\-F^\top\E[\check{y}_1]\-\widetilde F^\top\E[\check{\beta}_1^1]\-(A\+F)^\top \check y_2]dt,                                           \\
     & d\check{\varphi} \= -\{[A^{\!\top}\!\! \- \cPT\mathcal R^{\-1}B^\top]\check{\varphi} \+[F\-(\cPT\mathcal R^{\-1}D^\top \- C^\top)P\widetilde  F]\E[z]\-S\}dt, \\
     & S\=  Q\hat H\hat x \+ \hat H Q(\hat x\-\hat H\hat x)\-F^\top y_2\-F^\top\hat y_1 \-\widetilde F^\top\hat\beta_1,                                                                         \\
     & z_0\equiv 0,\quad \check\phi\equiv 0,\quad z_0(0)\=\xi_0,\quad z(0)\=\xi,\quad \check y_0(T)\=0,                                                                                      \\
     & \check{y}_1(T)\=0,\quad  \check{y}_2(T)\=0,\quad \check\phi(T)\=0,\quad
    \check{\varphi}(T)\=0,                                                                                                                                                                 \\
  \end{aligned}\right. \end{equation*}\end{shrinkeq}\n
which is consistent with the result of \cite[equation (33)]{QHX2020}. Moreover, $(\Lambda_1, \Lambda_2)$ takes the following form
\begin{shrinkeq}{-1ex}\begin{equation*} \left\{\begin{aligned}
     & \Lambda_1=-(R+D^\top\! PD)^{-1}(B^\top\! P+D^\top\! PC),                        \\
     & \Lambda_2=-(R+D^\top\! PD)^{-1}(B^\top\varphi+D^\top\! P\widetilde  F\hat x),\\
  \end{aligned}\right. \end{equation*}\end{shrinkeq}
which is also consistent with \cite{QHX2020}.

Through the discussion above, we compare our result with some previous literature. For the next part, we will study the performance of the mean-field {\g strategy}. Specifically, we will prove its  asymptotic optimality.
\section{Asymptotic $\varepsilon$-optimality}\label{asymptotic optimality}
\begin{definition}
  A mixed strategy set $\{ u_i^\varepsilon\in\mathcal U_i\}_{i=0}^N$ is called  asymptotically $\varepsilon$-optimal if there exists $\varepsilon=\varepsilon(N)>0$,
  $\lim_{N\rightarrow\infty}\varepsilon(N)=0$ such that
  \begin{shrinkeq}{-1ex}\begin{equation*}
    \left\{\begin{aligned}
      &\mathcal J_0( u_0^\varepsilon, u_{-0}^\varepsilon)\leq\inf_{u_0\in\mathcal U_0}\mathcal J_0(u_0, u_{-0}^\varepsilon)\+\varepsilon, \\
      &\scaleobj{1}{\frac{1}{N}}\big(\mathcal J_{soc}^{(N)}
      ( u_0^\varepsilon, u_{-0}^\varepsilon)-\inf_{u_{-0}\in\mathcal U_{-0}}\mathcal J_{soc}^{(N)}( u_0^\varepsilon,u_{-0})\big)\leq\varepsilon,
    \end{aligned}\right.
  \end{equation*}\end{shrinkeq}
  where $u_{-0}^\varepsilon:=\left\{u_{1}^\varepsilon, \cdots, u_{N}^\varepsilon\right\}$. In this case, $u_0^\varepsilon, u_{-0}^\varepsilon$ achieve an asymptotic $\varepsilon$-equilibrium, and $u_{1}^\varepsilon, \cdots, u_{N}^\varepsilon $ achieve an asymptotic $\varepsilon$-social optimum.
\end{definition}
Let $\widetilde u$ be the mean-field {\g strategy} given in Section V
and the realized decentralized states $(\widetilde x_0,\widetilde x_1,\cdots,\widetilde x_N)$ satisfy:
\begin{equation*}
\left\{\begin{aligned}
         & d\tilde{x}_0=(A_0\tilde{x}_0 \+ B_0\tilde{u}_0 \+ F_0\tilde{x}^{(N)})dt,\\
         & \hspace{1.1cm}+(C_0 \tilde{x}_0 \+ D_0\tilde{u}_0 \+ \widetilde F_0\tilde{x}^{(N)})dW_0, \quad \tilde{x}_0(0)=\xi_0                      \\
         & d\tilde{x}_i=(A\tilde{x}_i \+ B\tilde{u}_i \+ F\tilde{x}^{(N)})dt       \\
         & \hspace{1.1cm}+[C\tilde{x}_i+D\tilde{u}_i+\widetilde F\tilde{x}^{(N)} \+ \widetilde G\tilde{x}_0]dW_i, \quad \tilde{x}_i(0)=\xi.   \\
      \end{aligned}\right.\end{equation*}
where $\widetilde x^{(N)}=\frac{1}{N}\sum_{i=1}^N\widetilde x_i$.
  First, we need some estimations. In the proofs below, we will use $K$ to denote a generic constant whose value may change from line to line.
\begin{lemma}\label{estimation-1}\cite[Lemma 5.1]{HHN2018}
  Under (H1)-(H3) and (SA), there exists a constant $K_1$ independent of $N$ such that
  \begin{equation*}\rr
    \sup_{0\leq i\leq N}\E\sup_{0\leq t\leq T}\|\widetilde x_i(t)\|^2 \+ \sup_{0\leq i\leq N}\E\sup_{0\leq t\leq T}\|\z_i(t)\|^2\leq K_1.
  \end{equation*}
\end{lemma}
\begin{lemma}\label{estimation-2}
  Under (H1)-(H3) and (SA), there exists a constant $K_2$ independent of $N$ such that
  \begin{equation*}
    \begin{aligned}
      \E\sup_{0\leq t\leq T}\Big\|\widetilde x^{(N)}(t)-\Eo[z]\Big\|^2\leq \frac{K_2}{N}. \\
    \end{aligned}
  \end{equation*}
\end{lemma}
\begin{proof}
  It is easy to get that
  \begin{equation*}\begin{aligned}
        & d\Big(\widetilde x^{(N)}-\Eo[z]\Big)                                                                                              \\
      = & (A \+F)\Big(\widetilde x^{(N)}-\Eo[z]\Big) \+ B\Lambda_1\left(\frac{1}{N}\sum_{i=1}^{N}\z_i \- \Eo [z]\right) dt                                                                     \\
        & \+\frac{1}{N}\sum_{i=1}^N\Big[C\widetilde x_i\+D\tilde{u}_i\+\widetilde F\widetilde x^{(N)}\+\widetilde G\widetilde x_0\Big]dW_i.
    \end{aligned}\end{equation*}
  Therefore,
  \begin{equation*}\begin{aligned}
           & \E\sup_{0\leq s\leq t}\Big\|\widetilde x^{(N)}(s)-\Eo [z]\Big\|^2                                                            \\
      \leq & K\E\int_0^t\Big\|\widetilde x^{(N)}-\Eo[z]\Big\|^2 \+ \left\|\frac{1}{N}\sum_{i=1}^{N}\z_i \- \Eo [z]\right\|^2 ds\\
      &\rr\+ \frac{K}{N^2}\E\sup_{0\leq s\leq t}\Big\|\int_0^s\sum_{i=1}^N \big[C\widetilde x_i\+D\tilde{u}_i   \+\widetilde F\widetilde x^{(N)} \+\widetilde G\widetilde x_0\big]dW_i\big\|^2.                                      \\
    \end{aligned}\end{equation*}
    By\s
    \begin{equation*}
    \left\{\begin{aligned}
    & d \left(\frac{1}{N}\sum_{i=1}^{N}\z_i \- \Eo [z]\right) \= \left(\left(A\+B\Lambda_1\right)\left(\frac{1}{N}\sum_{i=1}^{N}\z_i \- \Eo [z]\right)\right)dt       \\
               &\hspace{1.7cm} \+\frac{1}{N}\sum_{i=1}^{N}(\left(C
                 \+D\Lambda_1\right)\z_i \+ D\Lambda_2\+\widetilde F\Eo[z]\+\widetilde G\z_0)dW_i, \\
    &\left(\frac{1}{N}\sum_{i=1}^{N}\z_i \- \Eo [z]\right)(0) \= 0.
    \end{aligned}\right.
    \end{equation*}
    it is easy to get
    \begin{equation}\label{211207_1}
    \begin{aligned}
    \E\sup_{0\leq s\leq t}\Big\|\frac{1}{N}\sum_{i=1}^{N}\z_i \- \Eo [z]\Big\|^2 = O\left(\frac{1}{N}\right).
    \end{aligned}
    \end{equation}
  Then by Burkholder-Davis-Gundy inequality, we have \f
  \begin{equation*}
     \begin{aligned}
           & \E\sup_{0\leq s\leq t}\Big\|\widetilde x^{(N)}(s)-\Eo[z]\Big\|^2                                                                                        \\
      \leq & K\E\int_0^t\Big\|\widetilde x^{(N)}-\Eo[z]\Big\|^2ds                                                                                                    \\
           &\rr +\frac{K}{N^2}\E\int_0^t\sum_{i=1}^N\Big\|C\widetilde x_i+D\tu_i+\widetilde F\widetilde x^{(N)}+\widetilde G\widetilde x_0\Big\|^2ds                              \\
      \leq & K\E\int_0^t\Big\|\widetilde x^{(N)}\-\Eo[z]\Big\|^2ds \+\frac{K}{N}\Big(1\+\sup_{0\leq i\leq N}\E\sup_{0\leq t\leq T}\|\widetilde x_i(t)\|^2\Big). \\
    \end{aligned} \end{equation*}\n
  Finally, it follows from Gronwall's inequality, and {Lemma \ref{estimation-1}} that there exists a constant $K_2$ independent of $N$ such that
  \begin{equation*}
    \begin{aligned}
      \E\sup_{0\leq t\leq T}\Big\|\widetilde x^{(N)}(t)-\Eo[z]\Big\|^2\leq \frac{K_2}{N}. \\
    \end{aligned}
  \end{equation*}
\end{proof}

\begin{lemma}\label{estimation-3}
  Under (H1)-(H3)  and (SA), there exists a constant $K_3$ independent of $N$ such that
  \begin{equation*}
    \sup_{0\leq i\leq N}\E\sup_{0\leq t\leq T}\Big\|\widetilde x_i(t)-\z_i(t)\Big\|^2\leq \frac{K_3}{N}.
  \end{equation*}
\end{lemma}
\begin{proof}
  It is easy to check that\f
  \begin{equation*}
      \begin{aligned}
      d(\widetilde x_i\-\z_i)= & \Big[A(\widetilde x_i\-\z_i)\+F(\widetilde x^{(N)} \- \Eo[z])\Big]dt                                                                             \\
                                      & \+\!\Big[C(\widetilde x_i \- \z_i) \+ \widetilde F(\widetilde x^{(N)}\-\Eo[z]) \+\widetilde G(\widetilde x_0\-\z_0)\Big]dW_i, \\
    \end{aligned} \end{equation*}\n
  and
  \begin{equation*} \begin{aligned}
            & d(\widetilde x_0-\z_0)                                                                                                      \\
          = & \Big[A_0(\widetilde x_0-\z_0)+F_0(\widetilde x^{(N)}-\Eo[z])\Big]dt                \\
            & +\Big[ C_0 (\widetilde x_0-\z_0)+\widetilde F_0(\widetilde x^{(N)}-\Eo[z])\Big]dW_0. \\
        \end{aligned} \end{equation*}
  Therefore, it follows from Burkholder-Davis-Gundy inequality that
  \begin{equation*}\resizebox{\linewidth}{!}{$\begin{aligned}
               & \E\sup_{0\leq s\leq t}\|\widetilde x_i(s)-\z_i(s)\|^2                                                                                                                                                           \\
          \leq & K\E\int_0^t\Big[\|\widetilde x_i -\z_i \|^2+\|\widetilde x^{(N)} -\Eo[z]\|^2\Big]ds                                                                                                  \\
               & +\! 2\E\!\!\!\sup_{0\leq s\leq t}\!\Big\|\!\int_0^s\!\!\Big[C(\widetilde x_i \- \z_i)(r) \+ \widetilde F(\widetilde x^{(N)}\!(r) \!\- \E[z|\F_r^{W_0}])          \\
               & \hspace{24mm}+\widetilde G(\widetilde x_0-\z_0)(r)\Big]dW_i(r)\Big\|^2                                                                                                                                                  \\
          \leq & K\E\!\!\int_0^t\!\!\|\widetilde x_i \- \z_i \|^2ds \+ K\E\!\!\int_0^t\!\! \Big[\|\widetilde x^{(N)} \!\- \Eo[z]\|^2 \+ \|\widetilde x_0 -\z_0 \|^2\!\Big]ds, \\
        \end{aligned}$}\end{equation*}
  and
  \begin{equation*}\begin{aligned}
           & \E\sup_{0\leq s\leq t}\|\widetilde x_0(s) \-\z_0(s)\|^2                                                                               \\
      \leq & K\E\int_0^t\Big[\|\widetilde x_0\-\z_0\|^2+\|\widetilde x^{(N)}\-\Eo[z]\|^2\Big]ds                      \\
           & +2\E\sup_{0\leq s\leq t}\Big\|\int_0^s\Big[C_0(\widetilde x_0\-\z_0)(r)                                                 \\
           & \hspace{25mm}+\!\widetilde F_0(\widetilde x^{(N)}(r)\-\Eo[z])\!\Big]dW_0(r)\Big\|^2                                    \\
      \leq & K\E\int_0^t\|\widetilde x_0\-\z_0\|^2ds+K\E\int_0^t\Big\|\widetilde x^{(N)} \-\Eo[z]\Big\|^2ds.
    \end{aligned}\end{equation*}
  Therefore, it follows from Gronwall's inequality and {Lemma \ref{estimation-2}} that
  \begin{equation*}
    \sup_{0\leq i\leq N}\E\sup_{0\leq t\leq T}\|\widetilde x_i(t)-\z_i(t)\|^2\leq \frac{K_3}{N}.\\
  \end{equation*}
\end{proof}
\subsection{Major agent}\label{major perturbation}
\begin{lemma}\label{lemma1}
  Under (H1)-(H3) and (SA),
  \begin{shrinkeq}{-1ex}\begin{equation*}
    |\mathcal J_0(\widetilde u_0,\widetilde u_{-0})-J_0(\bar v_0)|=O(\scaleobj{.7}{\frac{1}{\sqrt{N}}}).
  \end{equation*}\end{shrinkeq}
\end{lemma}
 \begin{proof}
  Recall \eqref{cost-major} and \eqref{limit cost-major}, it follows from\f
  \begin{shrinkeq}{-1ex}\begin{equation*} \begin{aligned}
           & \mathcal J_0(\widetilde u_0,\widetilde u_{-0})-J_0(\bar v_0)                                                                                                   \\
      =    & \frac{1}{2}\E \scaleobj{.8}{\int_0^T}\{\langle Q_0(\widetilde x_0-H_0\widetilde x^{(N)}),\widetilde x_0-H_0\widetilde x^{(N)}\rangle                       \\
           & \hspace{10mm}-\langle Q_0(\bar z_0-H_0\Eo [z]),\bar z_0-H_0\Eo [z]\rangle\}dt                  \\
      \leq & K\E \scaleobj{.8}{\int_0^T}[\|\widetilde x_0-\bar z_0\|^2+\|\widetilde x^{(N)}-\Eo [z]\|^2]dt                                        \\
           & \+K\scaleobj{.8}{\int_0^T}[(\E \|\widetilde x_0 \- \bar z_0\|^2)^\frac{1}{2} \+ (\E \|\widetilde x^{(N)} \- \Eo [z]\|^2)^\frac{1}{2}]dt  \=     O(\scaleobj{.7}{\frac{1}{\sqrt{N}}}),
    \end{aligned}\end{equation*}\end{shrinkeq}\n
where the last equality follows from Lemmas 1-3.
\end{proof}
{\g Consider the case when the major agent $\mathcal A_0$ uses an alternative strategy $u_0$ and the minor agent $\mathcal A_i$ still uses the strategy $\widetilde u_i$. The realized states with major agent's perturbation are}\f
\begin{shrinkeq}{-1ex}\begin{equation*}\left\{\label{state equation-major perturbation 1}\begin{aligned}
     & d\alpha_0\=(A_0 \alpha_0\+B_0u_0\+F_0  \alpha^{(N)})dt \+ (C_0 \alpha_0\+D_0u_0\+\widetilde F_0 \alpha^{(N)})dW_0, \\
     & d \alpha_i\=[A \alpha_i\+B\tilde{u}_i \+ F \alpha^{(N)}]dt \+[C\+D\tilde{u}_i
      \+\widetilde  F \alpha^{(N)}\+\widetilde G \alpha_0]dW_i,                       \\
     & \alpha_0(0)=\xi_0,\qquad \alpha_i(0)=\xi,\qquad 1\leq i\leq N,
  \end{aligned}\right.\end{equation*}\end{shrinkeq}\n
where $\alpha^{(N)}=\frac{1}{N}\sum_{i=1}^N\alpha_i$. The decentralized limiting states with major agent's perturbation are\s
\begin{shrinkeq}{-1ex}\begin{equation*}\left\{\label{state equation-major perturbation 2}\begin{aligned}
     & d\bar \alpha_0 \=(A_0\bar \alpha_0\+B_0 u_0\+F_0\Eo [z])dt\\
     & \hspace{10mm}\+(C_0\bar \alpha_0\+D_0u_0\+\widetilde F_0\Eo [z])dW_0,                                                                      \\
     & d\bar \alpha_i \=[A\bar \alpha_i\+B\tilde{u}_i\+F\Eo [z]]dt                                                                      \\
     & \hspace{10mm}\+[C\bar  \alpha_i \+ D\tilde{u}_i\+ \bar  F\Eo [z] \+ \widetilde G\bar \alpha_0]dW_i, \\
     & \bar \alpha_0(0)=\xi_0,\qquad {\g\bar \alpha_i}(0)=\xi,\qquad 1\leq i\leq N.
  \end{aligned}\right.\end{equation*}\end{shrinkeq}\n
\begin{lemma}\label{lemma2}
  Under (H1)-(H3) and (SA), we have
  \begin{shrinkeq}{-1ex}\begin{equation*}
    |\mathcal J_0( u_0,\widetilde u_{-0})-J_0( u_0)|=O(\scaleobj{.7}{\frac{1}{\sqrt{N}}}).
  \end{equation*}\end{shrinkeq}
\end{lemma}
\begin{proof}\f
  \begin{shrinkeq}{-1ex}\begin{equation*}  \begin{aligned}
           & \mathcal J_0( u_0,\widetilde u_{-0})-J_0( u_0)                                                                                                                                                  \\
      =    & \frac{1}{2}\E \scaleobj{.8}{\int_0^T}\Big\{\langle Q_0(\alpha_0-H_0\alpha^{(N)}),\alpha_0-H_0\alpha^{(N)}\rangle                                                                               \\
           & \hspace{12mm}-\langle Q_0(\bar\alpha_0-H_0\Eo [z]),\bar\alpha_0 -H_0\Eo [z]\rangle\Big\}dt                                         \\
      \leq & K\E \scaleobj{.8}{\int_0^T}\Big[\|\alpha_0-\bar\alpha_0\|^2 +\|\alpha^{(N)}-\Eo [z]\|^2\Big]dt                                                                                \\
           & +K\scaleobj{.8}{\int_0^T}\Big[(\E \|\alpha_0 \- \bar\alpha_0\|^2)^\frac{1}{2} \+ (\E \|\alpha^{(N)} \- \Eo [z]\|^2)^\frac{1}{2}\Big]dt \=    O(\scaleobj{.7}{\frac{1}{\sqrt{N}}}).\\
    \end{aligned} \end{equation*}\end{shrinkeq}\n
 \rr The last equality is similar to Lemmas 1-3.
\end{proof}
\begin{theorem}
Under (H1)-(H3) and (SA), $ \tilde u_0$ is an asymptotically $\varepsilon$-optimal strategy for the major agent.
\end{theorem}
\begin{proof}
  It follows from {Lemma \ref{lemma1}} and {Lemma \ref{lemma2}} that
  \begin{shrinkeq}{-1ex}\begin{equation*}
    \begin{aligned}
       & \mathcal J_0(\widetilde u_0,\widetilde u_{-0})\leq J_0(\bar v_0)+O(\scaleobj{.7}{\frac{1}{\sqrt{N}}})\leq J_0( u_0)+O(\scaleobj{.7}{\frac{1}{\sqrt{N}}}) \\
       & \hspace{18mm}\leq \mathcal J_0( u_0,\widetilde u_{-0})+O(\scaleobj{.7}{\frac{1}{\sqrt{N}}}).\\
    \end{aligned}
  \end{equation*}\end{shrinkeq}
\end{proof}
\subsection{Minor agents}\label{minor perturbation}
\subsubsection{Representation of social cost}\label{Representation of social cost}
Rewrite the large-population system \eqref{state equation-major} and \eqref{state equation-minor} as follows:
\begin{shrinkeq}{-1ex}\begin{equation}\label{LP system}
\resizebox{7.9cm}{!}{$
  \begin{aligned}
     & d\mathbf{x} \= (\mathbf{A}\mathbf{x} \+ \mathbf{B}u)dt \+\! \scaleobj{.8}{\sum_{i = 0}^{N}}(\mathbf{C}_i\mathbf{x} \+ \mathbf{D}_iu)dW_i,\quad \mathbf{x}(0) = {\Xi},
  \end{aligned}$}
\end{equation}\end{shrinkeq}
{\dg where\s
\begin{equation*}
  \begin{aligned}
     & \mathbf{A}=
    \left(\begin{smallmatrix}
        A_0 & \frac{F_0}{N} & \frac{F_0}{N}& \cdots &\frac{F_0}{N} \\
        0& A + \frac{F}{N} & \frac{F}{N} & \cdots & \frac{F}{N} \\
        0&\frac{F}{N} & A + \frac{F}{N} & \cdots & \frac{F}{N} \\
        \vdots & \vdots & \vdots & \ddots &\vdots\\
        0& \frac{F}{N} & \frac{F}{N} & \cdots & A + \frac{F}{N} \\
      \end{smallmatrix}\right),
    \mathbf{x}=
    \left(\begin{smallmatrix}
        x_0\\
        x_1 \\
        \vdots\\
        x_N
      \end{smallmatrix}\right), \\
     & \mathbf{B}=
    \left(\begin{smallmatrix}
        B_0& 0&0& \cdots & 0\\
        0&  B & 0 & \cdots & 0 \\
        0&  0 & B & \cdots & 0 \\
        \vdots& \vdots & \vdots & \ddots &\vdots\\
        0&  0 & 0 & \cdots & B \\
      \end{smallmatrix}\right),u=
    \left(\begin{smallmatrix}
        u_0\\
        u_1 \\
        \vdots\\
        u_N
      \end{smallmatrix}\right),\mathbf{C}_0= \left(\begin{smallmatrix}
        C_0& \frac{\widetilde F_0}{N}& \frac{\widetilde F_0}{N}& \cdots & \frac{\widetilde F_0}{N}\\
        0&  0 & 0 & \cdots & 0 \\
        0&  0 & 0 & \cdots & 0 \\
        \vdots& \vdots & \vdots & \ddots &\vdots\\
        0&  0 & 0 & \cdots & 0 \\
      \end{smallmatrix}\right) \\
     & \mathbf{D}_0= \left(\begin{smallmatrix}
        D_0& 0& 0& \cdots & 0\\
        0&  0 & 0 & \cdots & 0 \\
        0&  0 & 0 & \cdots & 0 \\
        \vdots& \vdots & \vdots & \ddots &\vdots\\
        0&  0 & 0 & \cdots & 0 \\
      \end{smallmatrix}\right),\mathbf{C}_i=\begin{smallmatrix}
      1 \\
      \vdots \\
      i+1 \\
      \vdots \\
      N+1
    \end{smallmatrix}
    \left(\begin{smallmatrix}
        0 & \cdots & 0 & 0 & 0 & \cdots & 0\\
        \vdots &  & \ddots & \vdots & \vdots & \vdots & \vdots\\
        \frac{\widetilde{F}}{N} & \cdots &\frac{\widetilde{F}}{N} & \frac{\widetilde{F}}{N} + C &\frac{\widetilde{F}}{N} &\cdots & \frac{\widetilde{F}}{N}\\
        \vdots & \vdots & \vdots & \vdots &  & \ddots & \vdots\\
        0 & \cdots & 0 & 0 & 0 & \cdots & 0\\
      \end{smallmatrix}\right),                                                          \\
     &
    \mathbf{D}_i=\begin{smallmatrix}
      1 \\
      \vdots \\
      i+1 \\
      \vdots \\
      N+1
    \end{smallmatrix}\left(\begin{smallmatrix}
        0 & \cdots & 0 & 0 & 0 &  \cdots & 0\\
        \vdots &  & \ddots & \vdots & \vdots & \vdots & \vdots\\
        0 & \cdots & 0  & D &0  & \cdots  & 0\\
        \vdots & \vdots & \vdots & \vdots &  & \ddots & \vdots\\
        0 & \cdots & 0  & 0 &0  & \cdots  & 0\\
      \end{smallmatrix}\right), \Xi = \left(\begin{smallmatrix}
        \xi_0 \\
        \xi\\
        \vdots\\
        \xi
      \end{smallmatrix}\right).
  \end{aligned}
\end{equation*}\n
Similarly, the social cost takes the following form:
\begin{equation}
  \begin{aligned}
     & \mathcal{J}^{(N)}_{soc}(u)                                                                                                         \\
     & = \frac{1}{2}\sum_{i=1}^n\E\int_0^T\Big[\Big\langle Q(x_i-Hx_0-\hat H x^{(N)}),                                            \\
     & \hspace{27mm}(x_i-Hx_0-\hat Hx^{(N)})\Big\rangle\+\langle Ru_i,u_i\rangle\Big]dt                                                \\
     & = \frac{1}{2}\E\int_{0}^{T}\Big[ \langle \mathbf{Q}\mathbf{x},\mathbf{x}\rangle \+\langle \mathbf{R}u,u\rangle \Big]dt, \\
  \end{aligned}
\end{equation}
where
\begin{equation*}
  \begin{aligned}
     & \mathbf{Q}\!=\!
    \left(\begin{smallmatrix}
        \mathbf{Q}_{00}& \mathbf{Q}_{01} & \mathbf{Q}_{02} &\cdots&\mathbf{Q}_{0N}\\
        \mathbf{Q}_{10}& \mathbf{Q}_{11} & \mathbf{Q}_{12} &  \cdots & \mathbf{Q}_{1N} \\
        \mathbf{Q}_{20}& \mathbf{Q}_{21} & \mathbf{Q}_{22} & \cdots & \mathbf{Q}_{2N} \\
        \vdots& \vdots & \vdots & \ddots &\vdots\\
        \mathbf{Q}_{N0}& \mathbf{Q}_{N1} & \mathbf{Q}_{N2} & \cdots & \mathbf{Q}_{NN} \\
      \end{smallmatrix}\right), \mathbf{R} \!=\! \left(\begin{smallmatrix}
        0&0&0&\cdots&0\\
        0& R & 0 & \cdots & 0\\
        0&0 & R & \cdots & 0\\
        \vdots&\vdots & \vdots & \ddots & \vdots\\
        0&0 & 0 & \cdots & R\\
      \end{smallmatrix}\right),
  \end{aligned}
\end{equation*}
and for $i=1,\cdots,N, j\neq i$,
\begin{equation*}
  \left\{
  \begin{aligned}
     & \mathbf{Q}_{00} = NQ+\hat H^\top Q\hat H - Q\hat H - \hat H^\top Q,                                                      \\
     & \mathbf{Q}_{0i} = -\hat H^\top QH + QH,                                                                                  \\
     & \mathbf{Q}_{i0} = -HQ\hat H+HQ,                                                                                          \\
     & \mathbf{Q}_{ii} = Q + \frac{1}{N}(\hat H^\top Q\hat H - Q\hat H - \hat H^\top Q), \mathbf{Q}_{ij} = \mathbf{Q}_{ii} - Q. \\
  \end{aligned}
  \right.
\end{equation*}}
Next,
define the following operators
\begin{shrinkeq}{-1ex}\begin{equation*}
  \resizebox{\linewidth}{!}{$
      \left\{
      \begin{aligned}
         & (Lu(\cdot))(\cdot) := \Phi(\cdot)\Big\{\scaleobj{.8}{\int_{0}^{\cdot}} \Phi(s)^{-1}\big[(\mathbf{B} - \scaleobj{.8}{\sum_{i = 0}^{N}}\mathbf{C}_i\mathbf{D}_i)u(s)\big]ds   \\
         & \hspace{3cm}+ \scaleobj{.8}{\sum_{i = 0}^{N} } \scaleobj{.8}{\int_{0}^{\cdot}}\Phi(s)^{-1} \mathbf{D}_iudW_i(s)\Big\},                                                       \\
         & \widetilde{L}u(\cdot) := (Lu(\cdot))(T),\quad \Gamma \Xi(\cdot) := \Phi(\cdot)\Phi^{-1}(0)\Xi,\quad \widetilde{\Gamma}\Xi := (\Gamma \Xi)(T).
      \end{aligned}
      \right.
    $}
\end{equation*}\end{shrinkeq}
Correspondingly, $L^*$ is defined as the adjoint operator of $L$.
Hence, we can rewrite the cost functional as follows:
\begin{shrinkeq}{-1ex}
\begin{align}
  \notag & 2\mathcal{J}^{(N)}_{soc}(u)
  =  \langle (L^*\mathbf QL \+ \mathbf{R})u(\cdot),u(\cdot)\rangle \+ 2\langle L^*Q\Gamma \Xi(\cdot) ,u(\cdot)\rangle            \\
  \notag & \hspace{20mm} \+\langle Q\Gamma \Xi(\cdot),\Gamma \Xi(\cdot)\rangle                                                                                 \\
         & \hspace{14.5mm} :=  \langle M_2u(\cdot),u(\cdot)\rangle \+ 2\langle M_1,u(\cdot)\rangle  \+ M_0 .
\end{align}
\end{shrinkeq}
Note that, $M_2$ is a self-adjoint positive semidefinite bounded linear operator.
{\dg
\subsubsection{Minor agent's perturbation}\label{Minor agent's perturbation}

Let us consider the case that the minor agent $\mathcal A_i$ uses an alternative strategy $u_i$ while the major agent and all other minor agents $\mathcal A_j,j\neq i$ use the strategies $\widetilde u_{-i}$. The realized states with the $i^\text{th}$ minor agent's perturbation are
\begin{equation*}\left\{\label{state equation-major perturbation 3}\begin{aligned}
     & d\check x_0=\Big[A_0  \check x_0\+B_0\tilde{u}_0\+F_0  \check x^{(N)}\Big]dt                    \\
     & \hspace{11mm}\+\Big[C_0\check x_0\+D_0\tilde{u}_0\+\widetilde F_0 \check x^{(N)}\Big]dW_0, \\
     & d \check x_i=\Big[A \check x_i\+B u_i\+F \check x^{(N)}\Big]dt                                                \\
     & \hspace{11mm}\+\Big[C  \check x_i\+D u_i
      \+\widetilde  F \check x^{(N)}\+\widetilde G \check x_0\Big]dW_i,                                              \\
     & d \check x_j=\Big[A  \check x_j\+B \tilde{u}_j\+F \check x^{(N)}\Big]dt                       \\
     & \hspace{11mm}\+\Big[C \check x_j \+ D  \tilde{u}_j
      \+\widetilde  F \check x^{(N)}\+\widetilde G \check x_0\Big]dW_j,                                              \\
     & y_0(0)=\xi_0,\qquad y_i(0)=\xi,\qquad 1\leq j\leq N,j\neq i,
  \end{aligned}\right.\end{equation*}
where $\check x^{(N)}=\frac{1}{N}\sum_{i=1}^N\check x_i$.
Moreover, introduce the following system:
\begin{equation*} \left\{\begin{aligned}
         & d\widetilde l_0=\Big[A_0\widetilde l_0+B_0\tu_0 \+ F_0\Eo[z]\Big]dt \\
         & \hspace{11mm}+\Big[C_0
        \widetilde l_0+D_0\tu_0+\widetilde F_0\Eo[z]\Big]dW_0,                            \\
         & d\widetilde l_j=\Big[A \widetilde l_j+B\tu_j+F\Eo[z]\Big]dt       \\
         & \hspace{11mm}+\Big[C
        \widetilde l_j+D\tu_j+\widetilde F\Eo[z]+\widetilde Gz_0\Big]dW_j,               \\
         & \widetilde l_0(0)=\xi_0,\qquad\widetilde l_j(0)=\xi,\qquad1\leq j\leq N,
      \end{aligned}\right. \end{equation*}
and\s
\begin{equation*} \left\{\label{state equation-major perturbation 4}\begin{aligned}
     & d\check l_0=\Big[A_0  \check l_0 +B_0 \tu_0 +F_0 \Eo[z]\Big]dt                 \\
     & \hspace{10mm}+\Big[ C_0  \check l_0+D_0\tu_0 +\widetilde F_0 \Eo[z]\Big]dW_0, \\
     & d \check l_i=\Big[A \check l_i+B u_i+F \check l^{(N)}\Big]dt                                                           \\
     & \hspace{10mm}+\Big[C  \check l_i+D u_i
      +\widetilde  F \check l^{(N)}+\widetilde Gz_0\Big]dW_i,                                                                 \\
     & d \check l_j=\Big[ A \check l_j+B \tu_j +F \Eo[z]\Big]dt                    \\
     & \hspace{10mm}+\Big[ C  \check l_j+ \tu_j
    +\widetilde  F \Eo[z]+\widetilde G z_0\Big]dW_j,                                               \\
     & \check l_0(0)=\xi_0,\quad \check l_i(0)=\check l_j(0)=\xi,\quad 1\leq j\leq N,\quad j\neq i,
  \end{aligned}\right.\end{equation*}\n
where $\widetilde l^{(N)}=\sum_{i=1}^N\frac{1}{N}\widetilde l_i$ and $\check l^{(N)}=\frac{1}{N}\sum_{i=1}^N\check l_i$.
Similar to {Lemma \ref{estimation-1}}, we have\f
\begin{equation*}\rr
  \left\{\begin{aligned}
      &  \E\sup_{0\leq t\leq T}\|\check x_i(t)\|^2   \+ \E\sup_{0\leq t\leq T}\|\check l_i(t)\|^2\leq K (1 \+ \|u_i\|^2_{L^2}),\\
      &\E\sup_{0\leq t\leq T}\|\widetilde l_i(t)\|^2 \leq K.
  \end{aligned}\right.
\end{equation*}\f
Similar to the computation in Section \ref{p-b-p optimality}, we have\s
\begin{equation*} \begin{aligned}
    \delta \mathcal J_{soc}^{(N)}= & \E\int_0^T\Big[\langle Q\widetilde x_i,\delta x_i\rangle+\langle R\widetilde u_i,\delta u_i\rangle-\Big\langle Q(\hat H\hat x+Hz_0)                                                                          \\
                                   & \hspace{6mm}+\!\hat H Q(\hat x \- \hat H\hat x \- Hz_0) \- F\!^\top\! y_2 \- F\!^\top\!\Eo[y_1^j]                                                                             \\
                                   & \hspace{6mm}-\!\widetilde F^\top\!\Eo[\beta_1^{jj}]\-F_0\!\!\!^\top\! y_1^0\-\widetilde F_0^\top\!\!\beta_1^0,\delta x_i\!\Big\rangle\Big] dt\+\!\sum_{l=1}^{13}\!\varepsilon_{l,i},
  \end{aligned}\end{equation*}\n
where $\varepsilon_{1,i},\cdots,\varepsilon_{13,i}$ are defined by \eqref{e1-4}, \eqref{e5-7}, \eqref{e8-10} and \eqref{e11-13}.
Finally, we have
\begin{equation*}\resizebox{\linewidth}{!}{$\begin{aligned}
        \delta \mathcal J_{soc}^{(N)}
        =\E\int_0^T & \Big[\langle Q\widetilde l_i,\delta l_i\rangle-\langle S,\delta l_i\rangle+\langle R\widetilde u_i,\delta u_i\rangle \Big] dt+\sum_{l=1}^{15}\varepsilon_{l,i}, \\
      \end{aligned}$}\end{equation*}
where\s
\begin{equation*}
   \left\{\begin{aligned}
     & \varepsilon_{14}=\E\int_0^T\Big[\langle Q(\widetilde x_i\-\widetilde l_i),\delta x_i\rangle\+\langle Q\widetilde l_i,\check x_i\-\check l_i\rangle \-\langle Q\widetilde l_i,\widetilde x_i\-\widetilde l_i\rangle\Big]dt, \\
     & \varepsilon_{15}=\E\int_0^T\Big[\langle S,\check x_i\-\check l_i\rangle\-\langle S,\widetilde x_i\-\widetilde l_i\rangle\Big] dt.\\
  \end{aligned}\right.\end{equation*}\n
At the end of this subsection, we give some lemmas which will be used in  subsection 3)
\begin{lemma}\label{estimation-8}
  Under (H1)-(H3) and (SA), there exists a constant $K_9$ independent of $N$ such that
  \begin{equation*}
    \E\sup_{0\leq t\leq T}\|\check x^{(N)}(t)-\Eo[z]\|^2\leq  K_9\left(\frac{1}{N} \+\frac{\| u_i\|^2_{L^2}}{N^2}  \right).
  \end{equation*}
\end{lemma}
\begin{proof}
  Note that\s
  \begin{equation*} \begin{aligned}
      d\check x^{(N)}= & \Big[A\check x^{(N)} + \frac{1}{N} B\Lambda_1\sum_{j=1}^{N}\z_j - \frac{1}{N}B\Lambda_1\z_i  + \frac{1}{N}Bu_i+\frac{N-1}{N}B\Lambda_2                             \\
                       & + F\check x^{(N)}\Big]dt+\frac{1}{N}\Big[C  \check x_i+D u_i
        +\widetilde  F \check x^{(N)}+\widetilde G \check x_0\Big]dW_i                                                                                             \\
                       & +\frac{1}{N}\sum_{j\neq i}\Big[C  \check x_j+D \tu_j + \widetilde  F \check x^{(N)}+\widetilde G \check x_0\Big]dW_j  \\
      =                & \Big[(A+F)\check x^{(N)} + \frac{1}{N} B\Lambda_1\sum_{j=1}^{N}\z_j  +\frac{1}{N}B(u_i-\Lambda_1\z_i-\Lambda_2)                                                       \\
                       & +B\Lambda_2\Big]dt+\frac{1}{N}\Big[C  \check x_i+D u_i
        +\widetilde  F \check x^{(N)}+\widetilde G \check x_0\Big]dW_i                                                                                             \\
                       & +\frac{1}{N}\sum_{j\neq i}\Big[C  \check x_j+D \tu_j + \widetilde  F \check x^{(N)}+\widetilde G \check x_0\Big]dW_j,\\
    \end{aligned}\end{equation*}\n
  and similar to \cite[equation (5.7)]{2} we have\s
  \begin{equation*}
    \begin{aligned}
       & d\Eo[z] =\Big[\Big(A-B\mathcal R^{-1}(B^\top P+D^\top PC)\Big)\Eo[z]                                                                                 \\
       & \hspace{14mm}-B\mathcal R^{-1}D^\top P\widetilde G\Eo[z_0]-\!B\mathcal R^{-1}\! B^\top\!\Eo[\check\varphi]                                           \\
       & \hspace{14mm}\-B\mathcal R^{-1}\! D^\top\!\!\Eo[\check\eta]\+\!\Big(\!F\-B\mathcal R^{-1}\!D^\top\!\! P\widetilde F\!\Big)\Eo[z]\!\Big]dt\\
    \end{aligned}
  \end{equation*}\n
  Thus,
  \begin{equation*}\begin{aligned}
       & d(\check x^{(N)}\-\Eo[z])                                                                                     \\=&\Big[(A+F)(\check x^{(N)}\-\Eo[z]) \+ B\Lambda_1\Big(\frac{1}{N}\sum_{j=1}^{N}\z_j \- \Eo(z)\Big)\\
       & +\frac{1}{N}B(u_i\-\Lambda_1 \z_i\-\Lambda_2)\Big]dt                                                                           \\
       & +\frac{1}{N}\Big[C  \check x_i+D u_i
        +\widetilde  F \check x^{(N)}+\widetilde G \check x_0\Big]dW_i                                                                             \\
       & +\frac{1}{N}\sum_{j\neq i}\Big[C  \check x_j+D \tu_j + \widetilde  F \check x^{(N)}+\widetilde G \check x_0\Big]dW_j. \\
    \end{aligned}\end{equation*}
  By Burkholder-Davis-Gundy inequality and \eqref{211207_1}, we have\f
  \begin{equation*}  \begin{aligned}
           & \E\sup_{0\leq s\leq t}\Big\|\check x^{(N)}(s)-\Eo[z]\Big\|^2                                                                                                                         \\\leq
           & K\E\int_0^t\Big[\Big\|\check x^{(N)} \- \Eo[z]\Big\|^2\!\!\+ \frac{1}{N^2} \|\check x_j\|^2 \!\+\frac{1}{N^2}\| u_i\|^2  \!\+ \frac{1}{N}\Big]ds\+O\Big(\frac{1}{N}\Big)  \\
           & \+\frac{K}{N^2}\E\sup_{0\leq s\leq t}\Big\|\int_0^s\Big[C  \check x_i\+D u_i \+\widetilde  F \check x^{(N)}\+\widetilde G \check x_0\Big]dW_i\Big\|^2                                               \\
           &\rr \+\frac{K}{N^2}\E\sup_{0\leq s\leq t}\Big\|\!\sum_{j\neq i} \!\int_{0}^{s}\! \Big[C  \check x_j\+D \tu_j
        \+\widetilde  F \check x^{(N)} \+\widetilde G \check x_0\Big]dW_j\Big\|^2                                                                                                                                              \\
      \leq & K\E\int_0^t\|\check x^{(N)}-\Eo[z]\|^2ds \+O\Big(\frac{1}{N}\Big) \+\frac{\| u_i\|^2_{L^2}}{N^2} .                                                                                                                \\
    \end{aligned}\end{equation*}\n
  Finally, it follows from Gronwall's inequality that
  \vspace{-2mm}
  \begin{equation*}
     \E\sup_{0\leq t\leq T}\|\check x^{(N)}(t)-\Eo[z]\|^2\leq K_9\left(\frac{1}{N} \+\frac{\| u_i\|^2_{L^2}}{N^2}  \right).\\
  \end{equation*}
\end{proof}
By {Lemma \ref{estimation-2}} and {Lemma \ref{estimation-8}}, we have
\begin{lemma}\label{estimation-9}
  Under (H1)-(H3) and (SA), there exist constants $K_{10}$, $K_{11}$, $ K_{12}$ independent of $N$ such that
  \vspace{-2mm}
  \begin{align}
     & \notag\E\sup_{0\leq t\leq T}\|\check x_0(t)-z_0(t)\|^2\leq  K_{10}\left(\frac{1}{N} \+\frac{\| u_i\|^2_{L^2}}{N^2}  \right),                                      \\
     &\notag \sup_{0\leq j\leq N}\E\sup_{0\leq t\leq T}\|\check x_j(t)-\check l_j(t)\|^2\leq  K_{11}\left(\frac{1}{N} \+\frac{\| u_i\|^2_{L^2}}{N^2}  \right),           \\
     &\notag \sup_{0\leq j\leq N}\E\sup_{0\leq t\leq T}\|\widetilde x_j(t)-\widetilde l_j(t)\|^2\leq K_{12}\left(\frac{1}{N} \+\frac{\| u_i\|^2_{L^2}}{N^2}  \right).
  \end{align}
\end{lemma}
\begin{lemma}\label{estimation-12}
  Under (H1)-(H3) and (SA), there exists a constant $K_{13}$ independent of $N$ such that
  \vspace{-2mm}
  \begin{equation*}
    \sup_{0\leq j\leq N,j\neq i}\E\sup_{0\leq t\leq T}\|\delta x_j\|^2\leq  K_{13}\left(\frac{1}{N^2} \+\frac{\| u_i\|^2_{L^2}}{N^2}  \right).
  \end{equation*}
\end{lemma}
\begin{proof}
  By the equations of $\delta x_0$, $\delta x_i$, $\delta x_j$ and Gronwall's inequality, we have
  \begin{align*}
     & \E\sup_{0\leq s\leq t}\|\delta x_0(s)\|^2\leq K\E\int_0^t\|\delta x^{(N)}\|^2ds,                             \\
     & \sup_{1\leq j\leq N,j\neq i}\E\sup_{0\leq s\leq t}\|\delta x_j(s)\|^2\leq K\E\int_0^t\|\delta x^{(N)}\|^2ds, \\
     & \E\sup_{0\leq s\leq t}\|\delta x_i(s)\|^2\leq K+K\E\int_0^t\|\delta x^{(N)}\|^2ds.                           \\[-6mm]
  \end{align*}
  It is easy to check that
  \begin{equation*}
    \resizebox{\linewidth}{!}{$
        \begin{aligned}
           & d\Big(\frac{1}{N}\sum_{j=1}^N\delta x_j\Big)=\Big[A\Big(\frac{1}{N}\sum_{j=1}^N\delta x_j\Big)+\frac{1}{N}B\delta u_i+F\delta x^{(N)}\Big]dt \\
           & \hspace{28mm}+\frac{1}{N}\sum_{j=1}^N\Big[C\delta x_j+\widetilde F\delta x^{(N)}+\widetilde G\delta x_0\Big]dW_j.\\
        \end{aligned}
      $}
  \end{equation*}
  Applying It\^{o}'s formula to $\|\sum_{j=1}^N\delta x_j\|^2$, we have\s
  \begin{equation*} \begin{aligned}
      \E\Big\|\frac{1}{N}\sum_{j=1}^N\delta x_j\Big\|^2 \!\leq\! & K\E \int_0^t \Big\|\frac{1}{N}\sum_{j=1}^N\delta x_j\Big\|^2+\frac{1}{N^2}\|\delta u_i\|^2_{L^2}ds                    \\
                                                                         & +\frac{1}{N^2}\sum_{j=1}^N\E\int_0^t\Big[C\delta x_j+\widetilde F\delta x^{(N)}+\widetilde G\delta x_0\Big]^2ds \\
      \leq                                                               & K\E \int_0^t \Big\|\frac{1}{N}\sum_{j=1}^N\delta x_j\Big\|^2+ K\left(\frac{1}{N^2}\+\frac{\|\delta u_i\|^2_{L^2}}{N^2}\right).
    \end{aligned}\end{equation*}\n
  Moreover, we have
  \begin{equation*}
    \E\sup_{0\leq t\leq T}\|\delta x_0(t)\|^2\leq  K\left(\frac{1}{N^2}\+\frac{\|\delta u_i\|^2_{L^2}}{N^2}\right),
  \end{equation*}
  and
  \begin{equation*}
    \sup_{1\leq j\leq N,j\neq i}\E\sup_{0\leq t\leq T}\|\delta x_j(t)\|^2\leq  K\left(\frac{1}{N^2}\+\frac{\|\delta u_i\|^2_{L^2}}{N^2}\right).
  \end{equation*}
\end{proof}
\begin{lemma}\label{estimation-13}
  There exist constants $K_{14},K_{15},K_{16}$ independent of $N$ such that
  \begin{align}
     & \notag \E\sup_{0\leq t\leq T}\|x^{**}-\delta x_{-(0,i)}\|^2\leq  K_{14}\left(\frac{1}{N}\+\frac{\|\delta u_i\|^2_{L^2}}{N}\right), \\
     & \notag \E\sup_{0\leq t\leq T}\|x_0^{*}-N\delta x_0\|^2\leq  K_{15} \left(\frac{1}{N}\+\frac{\|\delta u_i\|^2_{L^2}}{N}\right),      \\
     & \notag \E\sup_{0\leq t\leq T}\|x_j^{*}-N\delta x_j\|^2\leq  K_{16}\left(\frac{1}{N}\+\frac{\|\delta u_i\|^2_{L^2}}{N}\right).
  \end{align}
\end{lemma}
\begin{proof}
  First, we have the following dynamics\s
  \begin{equation*} \left\{\begin{aligned}
        & d(x^{**}-\delta x_{-(0,i)})                                                                 \\
      = & \Big[(A+F)(x^{**}-\delta x_{-(0,i)})+\frac{1}{N}(F\delta x_i+\delta x_{-(0,i)})\Big]dt      \\
        & -\sum_{j\neq i}\Big[C\delta x_j+\widetilde F\delta x^{(N)}+\widetilde G\delta x_0\Big]dW_j, \\\\
        & d(x_0^*-N\delta x_0)                                                                        \\
      = & \Big[A_0(x_0^*-N\delta x_0)+F_0(x^{**}-\delta x_{-(0,i)})\Big]dt                            \\
        & +\Big[C_0(x_0^*-N\delta x_0)+\widetilde F_0(x^{**}-\delta x_{-(0,i)})\Big]dW_0,             \\\\
        & d(x_j^*-N\delta x_j)                                                                        \\
      = & \Big[A(x_j^*-N\delta x_j)+F(x^{**}-\delta x_{-(0,i)})\Big]dt                                \\
        & +\Big[C(x_j^*-N\delta x_j)+\widetilde F(x^{**}-\delta x_{-(0,i)})                           \\
        & \hspace{6mm}+\widetilde G(x_0^*-N\delta x_0)\Big]dW_j.
    \end{aligned}\right.\end{equation*}\n
  Therefore, it follows from Burkholder-Davis-Gundy inequality that
  \vspace{-1mm}\f
  \begin{equation*}
      \begin{aligned}
           & \E\sup_{0\leq s\leq t}\|x^{**}(s)-\delta x_{-(0,i)}(s)\|^2                                                                                                                    \\
      \leq & K\E\int_0^t\Big[\|x^{**}(s)-\delta x_{-(0,i)}(s)\|^2\+\frac{1}{N^2}\|F\delta x_i\+\delta x_{-(0,i)}\|^2\Big]ds                                                          \\
           & \+K\E\sup_{0\leq s\leq t}\Big\|\sum_{j\neq i}\int_0^s(\delta x_j\+\delta x^{(N)}\+\delta x_0)dW_j\Big\|^2                                                            \\
      \leq & K\E\int_0^t\|x^{**}(s)-\delta x_{-(0,i)}(s)\|^2ds                                                                                                                             \\
           &\hspace{-2mm} \+K\E\!\!\int_0^t\!\!\Big[\frac{1}{N^2}\|F\delta x_i\+\delta x_{-(0,i)}\|^2\+\!\sum_{j\neq i}\|\delta x_j\|^2\+\|\delta x^{(N)}\|^2 \+\|\delta x_0\|^2\Big]ds. \\
    \end{aligned} \end{equation*}\n
  It then follows from Gronwall's inequality and {Lemma \ref{estimation-12}} that
  \vspace{-2mm}
  \begin{equation*}
    \E\sup_{0\leq t\leq T}\|x^{**}(t)-\delta x_{-(0,i)}(t)\|^2\leq  K\left(\frac{1}{N}\+\frac{\|\delta u_i\|^2_{L^2}}{N}\right).\\
  \end{equation*}
  Similarly, we have the other two estimations.
\end{proof}}
\subsubsection{Asymptotic optimality}
{\b In order to prove asymptotic optimality for the minor agents, it
suffices to consider the perturbations $u_{-0}\in\mathcal U_{-0}$ such
that \f$\mathcal J_{soc}^{(N)}(\widetilde u_0,u_{-0})\leq\mathcal
J_{soc}^{(N)}(\widetilde u_0,\widetilde u_{-0})$\n. It is easy to check
that \f$\mathcal J_{soc}^{(N)}(\widetilde u_0,\widetilde u_{-0})\leq
KN$\n,
where $K$ is a constant independent of $N$. Therefore, in what follows,
we only consider the perturbations $u_{-0}\in\mathcal U_{-0}$ satisfying
$\sum_{j=1}^N\E \int_0^T\|u_j\|^2dt\leq KN$.
Let $\delta u_i=u_i-\widetilde u_i$,
and consider a perturbation  ${u} = \widetilde{u} + (0,\delta
u_1,\cdots,\delta u_N):=\widetilde{ u}+ \delta u$. Then by Section
\ref{Representation of social cost}, we have}
\begin{shrinkeq}{-0.5ex}\begin{equation}
\resizebox{7.9cm}{!}{$
  \begin{aligned}
       &2\mathcal{J}^{(N)}_{soc}(\widetilde{u} \+ \delta u) \=  \langle M_2(\widetilde{u} \+ \delta u),\widetilde{u} \+ \delta u\rangle \+ 2\langle M_1,\widetilde{u} \+ \delta u\rangle  \+ M_0                            \\
    &\hspace{1cm}=  2\mathcal{J}^{(N)}_{soc}(\widetilde{u} )  \+ 2 \langle M_2\widetilde{u} \+ M_1 , \delta u\rangle \+ \langle M_2 \delta u, \delta u\rangle, \\
  \end{aligned}$}
\end{equation}\end{shrinkeq}
where $\langle M_2\widetilde{u} + M_1 , \delta u_i\rangle$ is the Fr\'{e}chet differential of ${\mathcal J}^{(N)}_{soc}$ on $\widetilde{u}$ with variation $\delta u_i$. Therefore, in order to  prove asymptotic optimality for the minor agents, we only need to show that $ \langle M_2\widetilde{u} \+ M_1 , \delta u\rangle=o(N)$.
To this end, we introduce another assumption:
\begin{description}
  \item[(H4)] There exists constants $L_1,L_2>0$ independent of $N$ such that
      \begin{shrinkeq}{-0.5ex}
        \begin{align}
           & \label{convergence of BSDE-y}\E \scaleobj{.8}{\int_0^T}\Big\|\Eo [y_1^1]-\scaleobj{.8}{\frac{1}{N}\sum_{j\neq i}y_1^j}\Big\|^2 dt\leq \scaleobj{.8}{\frac{L_1}{N}},               \\
           & \label{convergence of BSDE-z}\E \scaleobj{.8}{\int_0^T}\Big\|\Eo [\beta_1^{11}]-\scaleobj{.8}{  \frac{1}{N}\sum_{j\neq i}\beta_1^{jj}}\Big\|^2 dt\leq \scaleobj{.8}{\frac{L_2}{N}}.
        \end{align}
        \end{shrinkeq}
\end{description}
\begin{theorem}
Under (H1)-(H4) and (SA), $(\tilde u_1,\cdots,\tilde u_N)$   is an asymptotically $\varepsilon$-optimal strategy set for the minor agents {\r whose individual cost functionals are given by \eqref{cost-minor} and the social cost is given by \eqref{cost-social optimal}}.
\end{theorem}
\begin{proof} \dg
  From Section \ref{Minor agent's perturbation}, we have
  \begin{equation*}
    \begin{aligned}
        &  \langle M_2\widetilde u+M_1,\delta u\rangle                                                                                                                                 \\
      = & \sum_{i=1}^{N}\E\int_0^T\Big[\langle Q\widetilde l_i,\delta l_i\rangle-\langle S,\delta l_i\rangle+\langle R\widetilde u_i,\delta u_i\rangle \Big] dt+\sum_{i=1}^{N}\sum_{l=1}^{15}\varepsilon_{l,i}.
    \end{aligned}
  \end{equation*}
  From the optimality of $\widetilde u$, we have
  $$\E\int_0^T\Big[\langle Q\widetilde l_i,\delta l_i\rangle-\langle S,\delta l_i\rangle+\langle R\widetilde u_i,\delta u_i\rangle \Big] dt=0.$$
  \vspace{-1mm}
  Moreover, it follows from {Lemma \ref{estimation-1}-\ref{estimation-13} } and (H4) that
  $$\sum_{i=1}^{N}\sum_{l=1}^{15}\varepsilon_{l,i} =
    O\Big(\sqrt{N}\Big).$$
  Therefore, $ \langle M_2\widetilde{u} \+ M_1 , \delta u\rangle=O\Big(\sqrt{N}\Big)$.
\end{proof}
\quad\\
\begin{remark}
  Note that\f
 \begin{shrinkeq}{-1ex}\begin{equation*}
     \begin{aligned}
        & d\Big(\frac{1}{N}\sum_{j\neq i}\widetilde x_j\-\Eo [\widetilde x_1]\Big)                                                                                                               \\
      = & \Big[A\Big(\frac{1}{N}\sum_{j\neq i}\widetilde x_j\-\Eo [\widetilde x_1]\Big)\+F\left(\widetilde x^{(N)} \- \Eo [\widetilde x^{(N)}]\right)\\
      &\+ B\Lambda_1\Big(\frac{1}{N}\sum_{j\neq i}\bar{z}_j \- \Eo[\z_1] \Big) \- \frac{1}{N}B\Lambda_2 \- \frac{1}{N}F\tilde{x}^{(N)}\Big]dt \\
        & \+\frac{1}{N}\sum_{j\neq i}\Big[C\widetilde x_j \+ D\tu_j\+\widetilde F\tilde{x}^{(N)}\+\widetilde G\widetilde x_0\Big]dW_j.
    \end{aligned}\end{equation*}\end{shrinkeq} \n
  Therefore, it follows from Burkholder-Davis-Gundy inequality and Gronwall's inequality that\s
  \begin{shrinkeq}{-1ex}\begin{equation*}\begin{aligned}
       & \E \sup_{0\leq s\leq t}\Big\|\frac{1}{N}\sum_{j\neq i}\widetilde x_j(s)-\Eo [\widetilde x_1]\Big\|^2
      \leq \scaleobj{.9}{\frac{K}{N}}.
    \end{aligned}\end{equation*}\end{shrinkeq}\n
  If $C=0$, applying It\^{o}'s formula to \f $\big\|\frac{1}{N}\sum_{j\neq i}y_1^j-\Eo [y_1^1]\big\|^2$\n, it is easy to check that \eqref{convergence of BSDE-y} in (H4) holds.
\end{remark}
\begin{remark}
  If the state has the following form\s
  \begin{shrinkeq}{-1ex}\begin{equation*}
     \left\{\begin{aligned}
       & dx_0\=(A_0x_0\+B_0u_0\+F_0x^{(N)})dt \+(C_0x_0\+D_0u_0\+\widetilde F_0x^{(N)})dW_0, \\
       & dx_i\=(Ax_i\+Bu_i\+Fx^{(N)}\+Gx_0)dt \+DdW_i,\\
       & x_0(0)=\xi_0,\ x_i(0)=\xi,\ 1\leq i\leq N.
    \end{aligned}\right.\end{equation*}\end{shrinkeq}\n
  Assumption (H4) is not needed to obtain the asymptotic optimality of the minor agents. However, if the state equations of the minor agents take the form \eqref{state equation-minor}, we need to suppose the assumption (H4) hold and we will continue to study this in the future work.
\end{remark}
 Last, by combining Theorems 2, 3 and recalling the Definition 1, we have the following result
\begin{theorem}\label{them3}
{\r The mean-field strategies $\tilde{u}_0$, $\tilde{u}_{-0}$ achieve an asymptotic
  $\varepsilon$-equilibrium between the major agent and the aggregation of minor agents, where $\tilde{u}_0 = \Theta_1 \tilde{x}_0 +
  \Theta_2$, $\tilde{u}_{-0} \= (\tilde{u}_1,\cdots,\tilde{u}_N)$ and $\tilde{u}_i = \Lambda_1 \tilde{x}_i + \Lambda_2$. Moreover $\tilde{u}_{1}, \cdots,
  \tilde{u}_{N} $ achieve an asymptotic $\varepsilon$-social
  optimum among the aggregation of minor agents. Thus, $(\tilde{u}_0, \tilde{u}_{-0})$ is  asymptotically $\varepsilon$-optimal.}
%
%
%
\end{theorem}
\begin{proof}
  By Theorem 2, we have
  \begin{shrinkeq}{-1.5ex}\begin{equation*}
    \begin{aligned}
      \mathcal J_0( \tilde{u}_0, \tilde{u}_{-0})\leq\inf_{u_0\in\mathcal U_0}\mathcal J_0(u_0, \tilde{u}_{-0}) + O\Big(\scaleobj{.8}{\frac{1}{\sqrt{N}}}\Big). \\
    \end{aligned}
  \end{equation*}\end{shrinkeq}
  Moreover, by Theorem 3, we have\s
  \begin{shrinkeq}{-1.2ex}\begin{equation*}
    \begin{aligned}
     \scaleobj{1.2}{\frac{1}{N}}\Big(\mathcal J_{soc}^{(N)}
      ( \tilde{u}_0, \tilde{u}_{-0})-\inf_{u_{-0}\in\mathcal U_{-0}}\mathcal J_{soc}^{(N)}( \tilde{u}_0, u_{-0})\Big)\leq O\Big(\scaleobj{.8}{\frac{1}{\sqrt{N}}}\Big). \\
    \end{aligned}
  \end{equation*}\end{shrinkeq}\n
  Thus, by Definition 1, Theorem \ref{them3} holds straightforwardly.
\end{proof}
\section{Numerical analysis}
This section presents some numerical example to illustrate our theoretical results.
Our example is motivated by an electric charging control problem in presence of distributed information network. Relevant literature include e.g., \cite{GTL2013} and \cite{MCH2013}.

Consider two competitive charging providers in a power-grid network for given city. One provider (namely, the major agent in our model) still retains the traditional charging scheme upon centralized information, whereby its charging strategy
is determined by a central controller.
On the other hand, another provider, taking account the well-recognized distributed datum, prefers to adopt some decentralized charging scheme, where the
strategy is determined by each distributed charging unit on grid-node, only upon their decentralized information \emph{on or around that node}. In this case, such distributed provider is actually formalized into a cooperative team wherein all its sub-units or nodes acts as the minor agents as in our model.

Moreover, under some mild conditions on (linear) demand-supply curve, above competitive problem may be fit into a linear-quadratic setup whenever a quadratic deviation or tracing criteria is applied, as in \cite{FWC2016} and \cite{graber2016}. Thus, we cook one example in our theoretical framework.

 We now specify such example in details, by randomly generate its coefficients: \f
$A_0=\left(\begin{smallmatrix}
      0.6423 & 0.2057 \\
      0.0287 & 0.7907
    \end{smallmatrix}\right)$,
$F_0=\left(\begin{smallmatrix}
      0.9225 & 0.3780 \\
      0.6933 & 0.8048
    \end{smallmatrix}\right)$,
$C_0=\left(\begin{smallmatrix}
      0.0125 & 0.4517 \\
      0.4720 & 0.1117
    \end{smallmatrix}\right)$,
$\widetilde F_0=\left(\begin{smallmatrix}
      0.8084 & 0.4284\\
      0.7032 & 0.1955
    \end{smallmatrix}\right)$,
$A=\left(\begin{smallmatrix}
      0.1023 & 0.2995\\
      0.1027 & 0.9415
    \end{smallmatrix}\right)$,
$F=\left(\begin{smallmatrix}
      0.0377 & 0.8910\\
      0.2866 & 0.1003
    \end{smallmatrix}\right)$,
$C=\left(\begin{smallmatrix}
      0.1641 & 0.0360\\
      0.3271 & 0.8063
    \end{smallmatrix}\right)$,
$\widetilde F =\left(\begin{smallmatrix}
      0.3751 & 0.6241\\
      0.4491 & 0.5093
    \end{smallmatrix}\right)$,
$\widetilde G(\cdot) = \left(\begin{smallmatrix}
      0.5018 & 0.7881\\
      0.6989 & 0.1633
    \end{smallmatrix}\right)$,
$B_0 = \left(\begin{smallmatrix}
      0.4514 & 0.2916\\
      0.4309 & 0.9989
    \end{smallmatrix}\right)$,
$D_0(\cdot)= \left(\begin{smallmatrix}
      0.4514 & 0.2916\\
      0.4309 & 0.9989
    \end{smallmatrix}\right)$,
$B = \left(\begin{smallmatrix}
      0.4389 & 0.4766\\
      0.2411 & 0.3539
    \end{smallmatrix}\right)$,
$D = \left(\begin{smallmatrix}
      0.8756 & 0.9451\\
      0.7493 & 0.8354
    \end{smallmatrix}\right)$,
$Q_0 = \left(\begin{smallmatrix}
      0.6210 & 0\\
      0 & 0.8691
    \end{smallmatrix}\right)$,
$H_0 = \left(\begin{smallmatrix}
      0.3250 & 0\\
      0 & 0.5957
    \end{smallmatrix}\right)$,
$Q = \left(\begin{smallmatrix}
      0.8701 & 0\\
      0 & 0.1925
    \end{smallmatrix}\right)$,
$H = \left(\begin{smallmatrix}
      0.3865 & 0\\
      0 & 0.2957
    \end{smallmatrix}\right)$,
$\hat H = \left(\begin{smallmatrix}
      0.7027 & 0\\
      0 & 0.0354
    \end{smallmatrix}\right)$,
$R_0 = \left(\begin{smallmatrix}
      0.7160 & 0\\
      0 & 0.5594
    \end{smallmatrix}\right)$,
$R = \left(\begin{smallmatrix}
      0.3885 & 0\\
      0 & 0.4182
    \end{smallmatrix}\right)$,
$Q_0 = \left(\begin{smallmatrix}
      0.6210 & 0\\
      0 & 0.8691
    \end{smallmatrix}\right)$,
$Q = \left(\begin{smallmatrix}
      0.8701 & 0\\
      0 & 0.1925
    \end{smallmatrix}\right)$,
$R_0 = \left(\begin{smallmatrix}
      0.7160 & 0\\
      0 & 0.5594
    \end{smallmatrix}\right)$,
$R = \left(\begin{smallmatrix}
      0.3885 & 0\\
      0 & 0.4182
    \end{smallmatrix}\right)$. \n
It is easy to see that such generated coefficients are constants and surely $L^\infty$ and Lipschitz continuous.  Thus, assumption (H1)-(H2) hold.

In the following simulation, we will calculate the feedback
form mean-field strategies and also the corresponding state
trajectories of major and minor agents. The convergence of
the population average will be also simulated.

Firstly, we solve \eqref{CC} by decentralizing method and decoupling method.
{\dg As we mentioned above, the CC system can be rewritten as \eqref{CC-equivalent form}. Then by letting $\mathbb{X}_1:=\E[X|\mathcal F_t^{W_0}] $, $\mathbb{Y}_1:=\E[Y|\mathcal F_t^{W_0}]$, $\mathbb{X}_2:=(X - \mathbb{X}_1)$,  $\mathbb{Y}_2:=(Y - \mathbb{Y}_1)$.
\begin{equation} \left\{\begin{aligned}
&d\mathbb{X}_1 \=  \Big[(\mathbb A_1\+\bar{\mathbb A}_1) \mathbb{X}_1 \+\mathbb B_1\mathbb{Y}_1\+\mathbb F_1\E[Z_1|\mathcal F_t^{W_0}]\Big]dt\\  &\+\Big[(\mathbb C_1^0\+\bar{\mathbb C}_1^0)\mathbb{X}_1\+\mathbb D_1^0 E\mathbb{Y}_1\+\mathbb F_1^0\E[Z_1|\mathcal F_t^{W_0}]\Big] dW_0,\\
&\\
&d\mathbb{X}_2= \Big[\mathbb A_1 \mathbb{X}_2 \+\mathbb B_1\mathbb{Y}_2\+\mathbb F_1(Z_1 \- \E[Z_1|\mathcal F_t^{W_0}])\Big]dt\\
& \+\Big[\mathbb C_1^0\mathbb{X}_2 \+ \mathbb D_1^0 \mathbb{Y}_2 \+ \mathbb F_1^0 (Z_1 \- \E[Z_1|\mathcal F_t^{W_0}])\Big] dW_0\\
& \+\Big[\mathbb C_1^1\mathbb{X}_2 \+ (C_1^1 \+ \bar{\mathbb C}_1^1)\mathbb{X}_1 \+ \mathbb D_1^1\mathbb{Y}_2 \+ \mathbb D_1^1\mathbb{Y}_1  \\
& \+ \mathbb F_1^1(Z_1 \- \E[Z_1|\mathcal F_t^{W_0}]) \+ \mathbb F_1^1\E[Z_1|\mathcal F_t^{W_0}] \Big]dW_1,\\
&\\
&d\mathbb{Y}_1= \Big[(\mathbb A_2\+\bar{\mathbb A}_2)\mathbb{X}_1\+(\mathbb B_2\+\bar{\mathbb B}_2)\mathbb{Y}_1 \+  \mathbb{C}_2\E[Z_1|\mathcal F_t^{W_0}]  \\
&\+ ( \tilde{\mathbb C}_2 \+\bar{\mathbb C}_2)\E[Z_2|\mathcal F_t^{W_0}]
\Big]dt \+\E[Z_1|\mathcal F_t^{W_0}]dW_0 ,\\
&\\
&d \mathbb{Y}_2 \= \Big[\mathbb A_2 \mathbb{X}_2 \+ \mathbb B_2\mathbb{Y}_2 \+ \mathbb C_2(Z_1 \- \E[Z_1|\mathcal F_t^{W_0}]) \\
&\+  \tilde{\mathbb C}_2(Z_2 \- \E[Z_2|\mathcal F_t^{W_0}]) \Big]dt \+ (Z_1 \- \E[Z_1|\mathcal F_t^{W_0}])dW_0 \+ Z_2dW_1,\\
&\\
&X(0)=(\xi_0^\top,\xi^\top)^\top,\quad Y(T)=(0^\top,0^\top,0^\top,0^\top,0^\top)^\top.
\end{aligned}\right. \end{equation}
Further, by letting $\E[Z_1|\mathcal F_t^{W_0}] = \mathbb{Z}_1^1$, $(Z_1 - \mathbb{Z}_1^1) = \mathbb{Z}_1^2$, we can rewrite it as \f
\begin{equation} \left\{\begin{aligned}
&d\binom{\mathbb{X}_1}{\mathbb{X}_2}\\
& = \bigg[
\left(\begin{smallmatrix}
  (\mathbb A_1\+\bar{\mathbb A}_1) & 0 \\
  0 & \mathbb A_1
\end{smallmatrix}\right)\binom{\mathbb{X}_1}{\mathbb{X}_2} \+ \left(\begin{smallmatrix}
  \mathbb B_1 & 0 \\
  0 & \mathbb B_1
\end{smallmatrix}\right)\binom{\mathbb{Y}_1}{\mathbb{Y}_2} \+ \left(\begin{smallmatrix}
  \mathbb F_1 & 0 \\
  0 & \mathbb F_1
\end{smallmatrix}\right)\binom{\mathbb{Z}_1^1}{\mathbb{Z}_1^2}\bigg]dt\\ & \+ \bigg[
\left(\begin{smallmatrix}
  (\mathbb C_1^0\+\bar{\mathbb C}_1^0) & 0 \\
  0 & \mathbb C_1^0
\end{smallmatrix}\right)\binom{\mathbb{X}_1}{\mathbb{X}_2} \+ \left(\begin{smallmatrix}
  \mathbb D_1^0 & 0 \\
  0 &  \mathbb D_1^0
\end{smallmatrix}\right)\binom{\mathbb{Y}_1}{\mathbb{Y}_2} \+ \left(\begin{smallmatrix}
  \mathbb F_1^0 & 0 \\
  0 & \mathbb F_1^0
\end{smallmatrix}\right)\binom{\mathbb{Z}_1^1}{\mathbb{Z}_1^2}\bigg]dW_0
\\
& \+ \bigg[
\left(\begin{smallmatrix}
  0 & 0 \\
   (C_1^1 \+ \bar{\mathbb C}_1^1) & \mathbb C_1^1
\end{smallmatrix}\right)\binom{\mathbb{X}_1}{\mathbb{X}_2} \+ \left(\begin{smallmatrix}
  0 & 0 \\
   \mathbb D_1^1 & \mathbb D_1^1\end{smallmatrix}\right)\binom{\mathbb{Y}_1}{\mathbb{Y}_2} \+ \left(\begin{smallmatrix}
  0 & 0 \\
  \mathbb F_1^1 & \mathbb F_1^1
\end{smallmatrix}\right)\binom{\mathbb{Z}_1^1}{\mathbb{Z}_1^2}\bigg]dW_1,
\\
&\\
&d\binom{\mathbb{Y}_1}{\mathbb{Y}_2} = \bigg[
\left(\begin{smallmatrix}
  (\mathbb A_2\+\bar{\mathbb A}_2) & 0 \\
  0 & \mathbb A_2
\end{smallmatrix}\right)\binom{\mathbb{X}_1}{\mathbb{X}_2} \+ \left(\begin{smallmatrix}
  (\mathbb B_2\+\bar{\mathbb B}_2) & 0 \\
  0 & \mathbb B_2
\end{smallmatrix}\right)\binom{\mathbb{Y}_1}{\mathbb{Y}_2} \+ \left(\begin{smallmatrix}
  \mathbb C_2 & 0 \\
  0 & \mathbb C_2
\end{smallmatrix}\right)\binom{\mathbb{Z}_1^1}{\mathbb{Z}_1^2} \\
& \+ \left(\begin{smallmatrix}
  0 & 0 \\
  0 & \tilde{\mathbb C}_2
\end{smallmatrix}\right)\binom{0}{{Z}_2}\+ \left(\begin{smallmatrix}
  0 & (\tilde{\mathbb C}_2\+\bar{\mathbb C}_2) \\
  0 & - \tilde{\mathbb C}_2
\end{smallmatrix}\right)\E\bigg[\binom{0}{{Z}_2}|\mathcal F_t^{W_0}\bigg] \bigg]dt \+ \left(\begin{smallmatrix}
  I & 0 \\
  0 & I
\end{smallmatrix}\right)\binom{\mathbb{Z}_1^1}{\mathbb{Z}_1^2} dW_0\\
& \+ \binom{0}{Z_2}dW_1,
\\
&\\
&X(0)=(\xi_0^\top,\xi^\top)^\top,\quad Y(T)=(0^\top,0^\top,0^\top,0^\top,0^\top)^\top.
\end{aligned}\right. \end{equation} \n
By letting \f
\begin{equation} \left\{\begin{aligned}
&\mathbf{A}_1 = \left(\begin{smallmatrix}
  (\mathbb A_1+\bar{\mathbb A}_1) & 0 \\
  0 & \mathbb A_1
\end{smallmatrix}\right),\  \mathbf{B}_1 = \left(\begin{smallmatrix}
  \mathbb B_1 & 0 \\
  0 & \mathbb B_1
\end{smallmatrix}\right),\ \mathbf{C}_1 = \left(\begin{smallmatrix}
  \mathbb F_1 & 0 \\
  0 & \mathbb F_1
\end{smallmatrix}\right),\\
&\mathbf{A}_2 = \left(\begin{smallmatrix}
  (\mathbb C_1^0+\bar{\mathbb C}_1^0) & 0 \\
  0 & \mathbb C_1^0
\end{smallmatrix}\right),\  \mathbf{B}_2 = \left(\begin{smallmatrix}
  \mathbb D_1^0 & 0 \\
  0 &  \mathbb D_1^0
\end{smallmatrix}\right),\ \mathbf{C}_2 = \left(\begin{smallmatrix}
  \mathbb F_1^0 & 0 \\
  0 & \mathbb F_1^0
\end{smallmatrix}\right),\ \\
&\mathbf{A}_3 = \left(\begin{smallmatrix}
  0 & 0 \\
   (C_1^1 + \bar{\mathbb C}_1^1) & \mathbb C_1^1
\end{smallmatrix}\right),\  \mathbf{B}_3 = \left(\begin{smallmatrix}
  0 & 0 \\
   \mathbb D_1^1 & \mathbb D_1^1\end{smallmatrix}\right),\ \mathbf{C}_3 = \left(\begin{smallmatrix}
  0 & 0 \\
  \mathbb F_1^1 & \mathbb F_1^1
\end{smallmatrix}\right),\ \\
&\mathbf{A}_4 = \left(\begin{smallmatrix}
  (\mathbb A_2+\bar{\mathbb A}_2) & 0 \\
  0 & \mathbb A_2
\end{smallmatrix}\right),\  \mathbf{B}_4 = \left(\begin{smallmatrix}
  (\mathbb B_2+\bar{\mathbb B}_2) & 0 \\
  0 & \mathbb B_2
\end{smallmatrix}\right),\ \mathbf{C}_4^1 = \left(\begin{smallmatrix}
  \mathbb C_2 & 0 \\
  0 & \mathbb C_2
\end{smallmatrix}\right), \\
&\mathbf{C}_4^2 = \left(\begin{smallmatrix}
  0 & 0 \\
  0 & \tilde{\mathbb C}_2
\end{smallmatrix}\right),\  \mathbf{C}_4^3 = \left(\begin{smallmatrix}
  (\mathbb B_2+\bar{\mathbb B}_2) & 0 \\
  0 & \mathbb B_2
\end{smallmatrix}\right). \\
\end{aligned}\right. \end{equation} \n
We have
\begin{equation}\label{44} \left\{\begin{aligned}
&d\binom{\mathbb{X}_1}{\mathbb{X}_2} = \bigg[
\mathbf{A}_1\binom{\mathbb{X}_1}{\mathbb{X}_2} + \mathbf{B}_1\binom{\mathbb{Y}_1}{\mathbb{Y}_2} + \mathbf{C}_1\binom{\mathbb{Z}_1^1}{\mathbb{Z}_1^2}\bigg]dt\\ &\hspace{20mm} + \bigg[
\mathbf{A}_2\binom{\mathbb{X}_1}{\mathbb{X}_2} + \mathbf{B}_2\binom{\mathbb{Y}_1}{\mathbb{Y}_2} + \mathbf{C}_2\binom{\mathbb{Z}_1^1}{\mathbb{Z}_1^2}\bigg]dW_0
\\
&\hspace{20mm} + \bigg[
\mathbf{A}_3\binom{\mathbb{X}_1}{\mathbb{X}_2} + \mathbf{B}_3\binom{\mathbb{Y}_1}{\mathbb{Y}_2} + \mathbf{C}_3\binom{\mathbb{Z}_1^1}{\mathbb{Z}_1^2}\bigg]dW_1,
\\
&\\
&d\binom{\mathbb{Y}_1}{\mathbb{Y}_2} = \bigg[
\mathbf{A}_4\binom{\mathbb{X}_1}{\mathbb{X}_2} + \mathbf{B}_4\binom{\mathbb{Y}_1}{\mathbb{Y}_2} + \mathbf{C}_4^1\binom{\mathbb{Z}_1^1}{\mathbb{Z}_1^2} + \mathbf{C}_4^2\binom{0}{Z_2}\\
& + \mathbf{C}_4^3\E\bigg[\binom{0}{{Z}_2}|\mathcal F_t^{W_0}\bigg]\bigg]dt  + \left(\begin{smallmatrix}
  I & 0 \\
  0 & I
\end{smallmatrix}\right)\binom{\mathbb{Z}_1^1}{\mathbb{Z}_1^2} dW_0 + \binom{0}{Z_2}dW_1,
\\
&\\
&X(0)=(\xi_0^\top,\xi^\top)^\top,\  Y(T)=(0^\top,0^\top,0^\top,0^\top,0^\top)^\top.
\end{aligned}\right. \end{equation}
We can introduce the following Riccati equation
\begin{equation*}
\left\{
  \begin{aligned}
    &\dot{\mathbf{K}} + \mathbf{K}\mathbf{A}_1 + \mathbf{K}\mathbf{B}_1\mathbf{K} - \mathbf{A}_4 - \mathbf{B}_4\mathbf{K}+ (\mathbf{K}\mathbf{C}_1 - \mathbf{C}_4)\\
    &\times (\mathbf{I} - \mathbf{K}\mathbf{C}_2)^{-1}(\mathbf{K}\mathbf{A}_2 + \mathbf{K}\mathbf{B}_2\mathbf{K}) = 0,\\
    &\mathbf{K}(T) = 0.\\
  \end{aligned}
  \right.
\end{equation*}
Then it is easy to verify that
\begin{equation*}
\left\{
  \begin{aligned}
  &\binom{\mathbb{Y}_1}{\mathbb{Y}_2} = \mathbf{K}\binom{\mathbb{X}_1}{\mathbb{X}_2},\\
    &\binom{\mathbb{Z}_1^1}{\mathbb{Z}_1^2} = (\mathbf{I} - \mathbf{K}\mathbf{C}_2)^{-1}(\mathbf{K} \mathbf{A}_2 + \mathbf{K}\mathbf{B}_2\mathbf{K})\binom{\mathbb{X}_1}{\mathbb{X}_2}, \\
    &\binom{0}{{Z}_2} = (\mathbf{K} \mathbf{A}_3 + \mathbf{K}\mathbf{B}_3\mathbf{K})\binom{\mathbb{X}_1}{\mathbb{X}_2} + \mathbf{K}\mathbf{C}_3\binom{\mathbb{Z}_1^1}{\mathbb{Z}_1^2}.
  \end{aligned}
  \right.
\end{equation*}
Then \eqref{44} becomes a decoupled BSDE
\begin{equation}  \left\{\begin{aligned}
&d\binom{\mathbb{X}_1}{\mathbb{X}_2} = \bigg[
\mathbf{A}_1  \+ \mathbf{B}_1\mathbf{K}  \+ \mathbf{C}_1(\mathbf{I} \- \mathbf{K}\mathbf{C}_2)^{\-1}(\mathbf{K} \mathbf{A}_2 \+ \mathbf{K}\mathbf{B}_2\mathbf{K})\bigg]\binom{\mathbb{X}_1}{\mathbb{X}_2}dt\\ &  \+ \bigg[
\mathbf{A}_2 \+ \mathbf{B}_2\mathbf{K} \+ \mathbf{C}_2(\mathbf{I} \- \mathbf{K}\mathbf{C}_2)^{\-1}(\mathbf{K} \mathbf{A}_2 \+ \mathbf{K}\mathbf{B}_2\mathbf{K})\bigg]\binom{\mathbb{X}_1}{\mathbb{X}_2}dW_0
\\
&  \+ \bigg[
\mathbf{A}_3 \+ \mathbf{B}_3\mathbf{K} \+ \mathbf{C}_3(\mathbf{I} \- \mathbf{K}\mathbf{C}_2)^{\-1}(\mathbf{K} \mathbf{A}_2 \+ \mathbf{K}\mathbf{B}_2\mathbf{K})\bigg]\binom{\mathbb{X}_1}{\mathbb{X}_2}dW_1,
\\
&\\
&X(0)=(\xi_0^\top,\xi^\top)^\top.
\end{aligned}\right. \end{equation}
}
Then by using Euler–Maruyama method, Milstein method and Runge–Kutta method,  ($\hat x$, $z_0$, $y_1^0$, $\beta_1^0$, $\hat y_1$, $\hat\beta_1$, $y_2$) can be obtained. Further, by \eqref{Theta12} and \eqref{Lambda1,2}, we can calculate $\Theta_1$, $\Theta_2$, $\Lambda_1$, $\Lambda_2$. Then, the realized states can be obtained by \eqref{state equation-major} and \eqref{state equation-minor}. The following graphs are the first coordinate of the realized states.
\begin{shrinkeq}{-1ex}
\[\resizebox{\linewidth}{!}{\includegraphics[width=6.8cm]{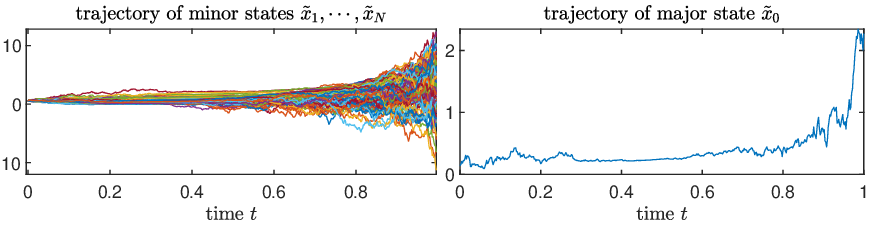}}\]
\end{shrinkeq}
Then, the corresponding feedback form mean-field controls can be obtained as well. The following graphs are the first coordinate of the mean-field controls.
\begin{shrinkeq}{-1ex}
\[\resizebox{\linewidth}{!}{\includegraphics[width=6.8cm]{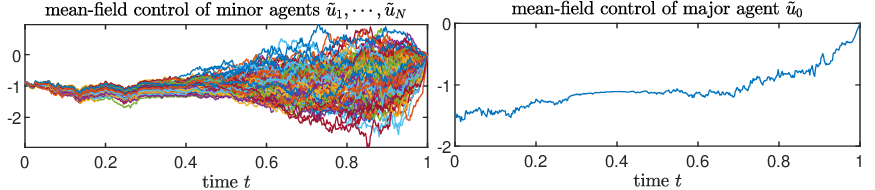}}\]
\end{shrinkeq}
Next, we simulate the convergence of the population state average $\check x^{(N)}(t)$ to the mean field $\hat{x}$. Specifically, we will calculate $\E \sup_{0\leq t\leq T}\|\check x^{(N)}(t)-\hat{x}\|^2$. First,  $\sup_{0\leq t\leq T}\|\check x^{(N)}(t)-\hat{x}\|^2$ can be calculated directly. Second, for the expectation, we repeat such process enough times (200 times) and take the average to simulate it.
\begin{shrinkeq}{-1ex}
\[\resizebox{\linewidth}{!}{\includegraphics[width=6.8cm]{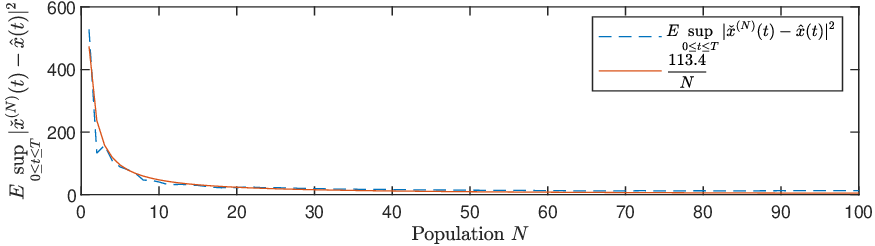}}\]
\end{shrinkeq}
The relation between \f$\E \sup_{0\leq t\leq T}\|\check x^{(N)}(t)-\hat{x}\|^2$\n and $N$ can be fitted by \f$\E \sup_{0\leq t\leq T}\|\check x^{(N)}(t)-\hat{x}\|^2 = \frac{113.4}{N}$\n with R-square 0.9944. In this sense, $\E \sup_{0\leq t\leq T}\|\check x^{(N)}(t)-\hat{x}\|^2 = O\big( \frac{1}{N}\big)$.

By the simulation above, we can see that the mean-field strategy is asymptotically optimal.

{\dg
\section{Appendix}

For any given $(Y,Z)$$\in$$ L_{\F}^2(0,T;\mathbb{R}^m)\!\times\! L_{\F}^2(0,T;\mathbb R^{m\times(d+1)})$ and $0\leq t\leq T$, the following SDE has a unique solution:
    \begin{equation}\label{1}
      \begin{aligned}
        X(t)= & x+\int_0^tb(s,X(s),\Eo[X(s)],Y(s),Z(s))ds                   \\
              & +\int_0^t\!\!\sigma(s,\!X(s),\Eo[X(s)],\! Y(s),Z(s))d W(s).
      \end{aligned}
    \end{equation}
    Therefore, we can introduce a map $\mathcal M_1: (Y,Z)\in L^2_{\F}(0,T;\mathbb R^m)\times L^2_{\F}(0,T;\mathbb R^{m\times(d+1)})\!\rightarrow\! X\in L^2_{\F}(0,T;\mathbb R^n)$ by \eqref{1}.
    Moreover, we have the following result:
    \begin{lemma}
      Let $X_i$ be the solution of \eqref{1} corresponding to $(Y_i,Z_i)$, $i=1,2$ respectively. Then for all $\rho\in\mathbb R$ and some constant $l_1>0$, we have\f
      \begin{equation}\label{2}
        \begin{aligned}
               & \E e^{-\rho t}\|\hat X(t)\|^2\+\bar\rho_1\E\!\!\int_0^te^{-\rho s}\|\hat X(s)\|^2ds                                                 \\
          \leq & (k_2l_1\+k_{11}^2)\E\!\!\int_0^te^{-\rho s}\|\hat Y(s)\|^2ds \+(k_3l_2\+k_{12}^2)\E\!\!\int_0^te^{-\rho s}\|\hat Z(s)\|^2ds,
        \end{aligned}
      \end{equation}\n and
      \begin{equation}\label{3}\resizebox{\linewidth}{!}{$
            \begin{aligned}
               & \E e^{-\rho t}\|\hat X(t)\|^2\leq(k_2l_1+k_{11}^2)\E\!\!\int_0^te^{-\bar\rho_1(t-s)-\rho s}\|\hat Y(s)\|^2ds \\
               & \hspace{25mm}+(k_2l_2+k_{12}^2)\E\!\!\int_0^te^{-\bar\rho_1(t-s)-\rho s}\|\hat Z(s)\|^2ds,
            \end{aligned}$}
      \end{equation}
      where $\bar\rho_1=\rho-2\rho_1-2k_1-k_2l_1^{-1}-k_3l_2^{-1}-k_9^2-k_{10}^2$ and $\hat\Phi=\Phi_1-\Phi_2$, $\Phi=X,Y,Z$. Moreover,
      \begin{equation}\label{4}
        \begin{aligned}
               & \E\int_0^Te^{-\rho t}\|\hat X(t)\|^2dt                                                            \\
          \leq & \frac{1-e^{-\bar\rho_1T}}{\bar\rho_1}\Big[(k_2l_1+k_{11}^2)\E\int_0^te^{-\rho t}\|\hat Y(t)\|^2dt \\
               & +(k_3l_2+k_{12}^2)\E\int_0^te^{-\rho t}\|\hat Z(t)\|^2dt\Big],
        \end{aligned}
      \end{equation}
      and
      \begin{equation}\label{5}
        \resizebox{\linewidth}{!}{$
            \begin{aligned}
               & e^{-\rho T}\E\|\hat X(T)\|^2\leq(1\vee e^{-\bar\rho_1T})\Big[(k_2l_1+k_{11}^2)\E\int_0^Te^{-\rho t}\|\hat Y(t)\|^2dt \\
               & \hspace{27mm}+(k_3l_2+k_{12}^2)\E\int_0^Te^{-\rho t}\|\hat Z(t)\|^2dt\Big].
            \end{aligned}
          $}
      \end{equation}
      Specifically, if $\bar\rho_1>0$,
      \begin{equation}\label{6}
        \begin{aligned}
           & e^{-\rho T}\E\|\hat X(T)\|^2\leq(k_2l_1+k_{11}^2)\E\int_0^Te^{-\rho t}\|\hat Y(t)\|^2dt \\
           & \hspace{27mm}+(k_3l_2+k_{12}^2)\E\int_0^Te^{-\rho t}\|\hat Z(t)\|^2dt.
        \end{aligned}
      \end{equation}
    \end{lemma}
    \begin{proof}
      For any $\rho>0$, applying It\^{o}'s formula to $e^{-\rho t}\|\hat X(t)\|^2$,\s
      \begin{equation*}\begin{aligned}
               & \E e^{-\rho t}\|\hat X(t)\|^2+\rho\E\int_0^te^{-\rho s}\|\hat X(s)\|^2ds                  \\
          =    & 2\E\int_0^te^{-\rho s}\hat X(s)(b(s,X_1(s),\Eo[X_1(s)],Y_1(s),Z_1(s)) \\
               & \hspace{10mm}-b(s,X_2(s),\Eo[X_2(s)],Y_2(s),Z_2(s)))ds                        \\
               & +\E\int_0^te^{-\rho s}(\sigma(s,X_1(s),\Eo[X_1(s)],Y_1(s),Z_1(s))     \\
               & \hspace{10mm}-\sigma(s,X_2(s),\Eo[X_2(s)],Y_2(s),Z_2(s)))^2ds                 \\
          \leq & \E\int_0^te^{-\rho s}\Big[(2\rho_1+2k_1+k_2l_1^{-1}+k_3l_2^{-1}+k_9^2+k_{10}^2)\|\hat X(s)\|^2   \\
               & \hspace{17mm}+(k_2l_1+k_{11}^2)\|\hat Y(s)\|^2+(k_3l_2+k_{12}^2)\|\hat Z(s)\|^2\Big]ds.
        \end{aligned}\end{equation*}\n
      Similarly, applying It\^{o}'s formula to $e^{-\bar\rho_1(t-s)-\rho s}\|\hat X(s)\|^2$, we have \eqref{3}.
      Integrating from $0$ to $T$ on both sides of \eqref{3} and noting that $\frac{1-e^{-\bar\rho_1(t-s)}}{\bar\rho_1}\leq \frac{1-e^{-\bar\rho_1T}}{\bar\rho_1},$
      we have
      \begin{equation*}\begin{aligned}
               & \E\int_0^Te^{-\rho t}\|\hat X(t)\|^2dt                                                         \\[-2mm]
          \leq & (k_2l_1+k_{11}^2)\frac{1-e^{-\bar\rho_1T}}{\bar\rho_1}\E\int_0^Te^{-\rho s}\|\hat Y(s)\|^2ds   \\
               & +(k_3l_2+k_{12}^2)\frac{1-e^{-\bar\rho_1T}}{\bar\rho_1}\E\int_0^Te^{-\rho s}\|\hat Z(s)\|^2ds.
        \end{aligned}\end{equation*}
      Letting $t=T$ in \eqref{3}, we have \eqref{5}.
    \end{proof}
    For any given $X\in L_{\F}^2(0,T;\mathbb R^n)$, the following BSDE has a unique solution:
    \begin{align}
         & \notag Y(t) = \int_t^T\!\!f(s, X(s), \Eo[X(s)], Y(s), \Eo[Y(s)], \\[-2mm]
         \label{7}& \hspace{18mm}Z(s), \Eo[Z(s)])ds - \int_t^T\!\!Z(s)dW(s).
    \end{align}
    Thus, we can introduce another map $\mathcal M_2\!:\!X\!\!\in\!\! L^2_{\F}(0,\!T;\!\mathbb R^n\!)\!\rightarrow \!(Y,Z)\in L^2_{\F}(0,T;\mathbb R^m)\!\times\! L^2_{\F}(0,\!T;\mathbb R^{m\!\times\!(d\+1)})$ by \eqref{7}.
    Similarly, we have the following result:
    \begin{lemma}
      Let $(Y_i,Z_i)$ be the solution of \eqref{7} corresponding to $X_i,i=1,2$, respectively. Then for all $\rho\in\mathbb R$ and some constants $l_3,l_4,l_5,l_6>0$ such that
      \begin{equation}\label{8}\begin{aligned}
               & \E e^{-\rho t}\|\hat Y(t)\|^2+\bar\rho_2\E\int_t^Te^{-\rho s}\|\hat Y(s)\|^2ds \\
               & +(1-k_7l_5-k_8l_6)\E\int_t^Te^{-\rho s}\|\hat Z(s)\|^2ds                               \\
          \leq & (k_4l_3+k_5l_4)\E\int_t^Te^{-\rho s}\|\hat X(s)\|^2ds,
        \end{aligned}\end{equation}
      and
      \begin{align}
        \notag         & \E e^{-\rho t}\|\hat Y(t)\|^2+(1-k_7l_5-k_8l_6)\E\int_t^Te^{-\rho s}\|\hat Z(s)\|^2ds \\
        \label{9} \leq & (k_4l_3+k_5l_4)\E\int_t^Te^{-\bar\rho_2(s-t)-\rho s}\|\hat X(s)\|^2ds,
      \end{align}
      where $\bar\rho_2=-\rho-2\rho_2-2k_6-k_4l_3^{-1}-k_5l_4^{-1}-k_7l_5^{-1}-k_8l_6^{-1}$,
      and $\hat\Phi=\Phi_1-\Phi_2$, $\Phi=X,Y,Z$. Moreover,\f
      \begin{equation}\label{10}
         \begin{aligned}
           & \E \!\! \int_0^Te^{-\rho t}\|\hat Y(t)\|^2dt \leq \frac{1-e^{-\bar\rho_2T}}{\bar\rho_2}(k_4l_3+k_5l_4)\E\!\!\int_0^Te^{-\rho s}\|\hat X(s)\|^2ds,
        \end{aligned}\end{equation}\n
      \vspace{-2mm}
      and
      \begin{equation}\label{11}\begin{aligned}
               & \E \int_0^Te^{-\rho t}\|\hat Z(t)\|^2dt                                                                                            \\
          \leq & \frac{(k_4l_3+k_5l_4)(1\vee e^{-\bar\rho_2T})}{(1-k_7l_5-k_8l_6)(1\wedge e^{-\bar\rho_2T})}\E\int_0^Te^{-\rho s}\|\hat X(s)\|^2ds.
        \end{aligned}\end{equation}
      Specifically, if $\bar\rho_2>0$,
      \begin{equation}\label{12}\resizebox{\linewidth}{!}{$\begin{aligned}
              \E \int_0^Te^{-\rho t}\|\hat Z(t)\|^2dt
              \leq\frac{k_4l_3+k_5l_4}{1-k_7l_5-k_8l_6}\E\int_0^Te^{-\rho s}\|\hat X(s)\|^2ds. \\[-1mm]
            \end{aligned}$}\end{equation}
    \end{lemma}
    \vspace{-2mm}
    \textbf{Proof of Theorem \ref{discounting}:}
    Define $\mathcal M:=\mathcal M_2\circ\mathcal M_1$, where $\mathcal M_1$ is defined by \eqref{1} and $\mathcal M_2$ is defined by \eqref{7}. Thus, $\mathcal M$ is a mapping from $L_{\F}^2(0,T;\mathbb R^m)\!\times\! L_{\F}^2(0,T;\mathbb R^{m\times (d+1)})$ into itself.
    For $(U_i,\!V_i)\!\in\! L_{\F}^2(0,\!T;\mathbb R^m)\!\times\! L_{\F}^2(0,\!T;\mathbb R^{m\!\times\! (d\+1)})$, let $X_i:=\mathcal M_1(U_i,V_i)$ and $(Y_i,Z_i):=\mathcal M(U_i,V_i)$.
    Therefore,
    \vspace{-2mm}
    \begin{equation*}\resizebox{\linewidth}{!}{$\begin{aligned}
                 & \E\int_0^T e^{-\rho t}\|Y_1(t)-Y_2(t)\|^2dt+
            \E\int_0^T e^{-\rho t}\|Z_1(t)-Z_2(t)\|^2dt                                                                  \\
            \leq & \Big[\frac{1-e^{-\bar\rho_2T}}{\bar\rho_2}+
            \frac{1\vee e^{-\bar\rho_2T}}{(1-k_7l_5-k_8l_6)(1\wedge e^{-\bar\rho_2T})}\Big](k_4l_3+k_5l_4)                       \\
                 & \times\E\int_0^T e^{-\rho t}\|X_1(t)-X_2(t)\|^2dt                                                     \\
            \leq & \Big[\frac{1-e^{-\bar\rho_2T}}{\bar\rho_2}+
            \frac{1\vee e^{-\bar\rho_2T}}{(1-k_7l_5-k_8l_6)(1\wedge e^{-\bar\rho_2T})}\Big]\frac{1-e^{-\bar\rho_1T}}{\bar\rho_1} \\
                 & \times(k_4l_3+k_5l_4)\Big[(k_2l_1+k_{11}^2)\E\int_0^T e^{-\rho t}\|U_1(t)-U_2(t)\|^2dt                \\[-2mm]
                 & +(k_3l_2+k_{12}^2)\E\int_0^T e^{-\rho t}\|V_1(t)-V_2(t)\|^2dt\Big].                                   \\[-2mm]
          \end{aligned}$}\end{equation*}
    Choosing suitable $\rho$, we get that $\mathcal M$ is a contraction mapping.
    \vspace{-1mm}

    Furthermore, if $2\rho_1+2\rho_2<-2k_1-2k_6-2k_7^2-2k_8^2-k_9^2-k_{10}^2$, we can choose $\rho\in\mathbb R$, $0<k_7l_5<\frac{1}{2}$ and $0<k_8l_6<\frac{1}{2}$ and sufficient large $l_1,l_2,l_3,l_4$ such that
\vspace{-2mm}
\begin{equation*}
  \begin{aligned}
    \bar\rho_1>0,\qquad \bar\rho_2>0,\qquad 1-k_7l_5-k_8l_6>0. \\[-2mm]
  \end{aligned}
\end{equation*}
Therefore,
\vspace{-2mm}
\begin{equation*}\resizebox{\linewidth}{!}{$\begin{aligned}
             & \E\int_0^T e^{-\rho t}\|Y_1(t)-Y_2(t)\|^2dt+
        \E\int_0^T e^{-\rho t}\|Z_1(t)-Z_2(t)\|^2dt                                            \\
        \leq & \Big[\frac{1}{\bar\rho_2}+
          \frac{1}{1-k_7l_5-k_8l_6}\Big]\frac{1}{\bar\rho_1}
        (k_4l_3+k_5l_4)                                                                                \\
             & \times\Big[(k_2l_1+k_{11}^2)\E\int_0^T e^{-\rho t}\|U_1(t)-U_2(t)\|^2dt         \\[-1mm]
             & \hspace{6mm}+(k_3l_2+k_{12}^2)\E\int_0^T e^{-\rho t}\|V_1(t)-V_2(t)\|^2dt\Big]. \\[-2mm]
      \end{aligned}$}\end{equation*}
The proof is complete. }


\begin{thebibliography}{00}


  \bibitem{B1996} F . Bloch,  “Sequential Formation of Coalitions in Games with Externalities and Fixed Payoff Division[J],” Games and Economic Behavior, vol. 14, no. 1, pp. 90–123, 1996.

  \bibitem{BSYS2016} A. Bensoussan, K. C. J. Sung, S. C. P. Yam, and S.-P. Yung, “Linear-quadratic mean field games,” J. Optim. Theory Appl., vol. 169, no. 2, pp. 496–529, 2016.



  \bibitem{CD2013} R. Carmona and F. Delarue, “Probabilistic Analysis of Mean-Field Games,” SIAM J. Control Optim., vol. 51, no. 4, pp. 2705–2734, Jan. 2013.

  \bibitem{carmona2017alternative}R. Carmona and P. Wang, “An alternative approach to mean field game with major and minor players, and applications to herders impacts,” Applied Mathematics \& Optimization, vol. 76, no. 1, pp. 5–27, 2017.

  \bibitem{carmona2016probabilistic} R. Carmona and X. Zhu, “A probabilistic approach to mean field games with major and minor players,” Ann. Appl. Probab., vol. 26, no. 3, pp. 1535–1580, Jun. 2016, Accessed: Aug. 21, 2020.




  \bibitem{FWC2016} H. Fang, Y. Wang, and J. Chen, “Health-aware and user-involved battery charging management for electric vehicles: Linear quadratic strategies,” IEEE Transactions on Control Systems Technology, vol. 25, no. 3, pp. 911-923, 2016.

  \bibitem{GTL2013} L. Gan, U. Topcu, and S. Low. “Optimal Decentralized Protocol for Electric Vehicle Charging," IEEE Trans. Power Sys., vol. 28, no, 2, pp. 940-951, 2013.

 \bibitem{graber2016} P.J. Graber,  “Linear Quadratic Mean Field Type Control and Mean Field Games with Common Noise, with Application to Production of an Exhaustible Resource," Appl. Math Optim., vol. 74, pp. 459–486, 2016.

  \bibitem{HK1983} S. Hart, and M. Kurz, “Endogenous formation of coalitions,” Econometrica: Journal of the econometric society, pp. 1047-1064, 1983.

  \bibitem{HHL2017} Y. Hu, J. Huang, and X. Li, “Linear quadratic mean field game with control input constraint,” ESAIM Control Optim. Calc. Var., vol. 24, no. 2, pp. 901–919, 2018, Accessed: Aug. 10, 2020.

  \bibitem{HHN2018} Y. Hu, J. Huang, and T. Nie, “Linear-Quadratic-Gaussian Mixed Mean-Field Games with Heterogeneous Input Constraints,” SIAM J. Control Optim., vol. 56, no. 4, pp. 2835–2877, Jan. 2018.

  \bibitem{HWY2019} J. Huang, B. Wang, and J. Yong, “Social Optima in Mean Field Linear-Quadratic-Gaussian Control with Volatility Uncertainty,” arXiv preprint arXiv:1912.06371, 2019.

  \bibitem{HWW2016} J. Huang, S. Wang, and Z. Wu, “Backward Mean-Field Linear-Quadratic-Gaussian (LQG) Games: Full and Partial Information,” IEEE Trans. Automat. Contr., vol. 61, no. 12, pp. 3784–3796, Dec. 2016.

  \bibitem{Huang2010} M. Huang, “Large-Population LQG Games Involving a Major Player: The Nash Certainty Equivalence Principle,” SIAM J. Control Optim., vol. 48, no. 5, pp. 3318–3353, Jan. 2010.

  \bibitem{HCM2007} M. Huang, P. E. Caines, and R. P. Malhame, “Large-Population Cost-Coupled LQG Problems With Nonuniform Agents: Individual-Mass Behavior and Decentralized $\varepsilon$-Nash Equilibria,” IEEE Trans. Automat. Contr., vol. 52, no. 9, pp. 1560–1571, Sep. 2007.

  \bibitem{HCM2012} M. Huang, P. E. Caines, and R. P. Malhame, “Social Optima in Mean Field LQG Control: Centralized and Decentralized Strategies,” IEEE Trans. Automat. Contr., vol. 57, no. 7, pp. 1736–1751, Jul. 2012.

  \bibitem{HCM2015} M. Huang, P. E. Caines, and R. P. Malhame, “Mean Field Games,” Handbook of Dynamic Game Theory. Springer, 2015.


  \bibitem{LL2007} J.-M. Lasry and P.-L. Lions, “Mean field games,” Jpn. J. Math., vol. 2, no. 1, pp. 229–260, 2007.

  \bibitem{MCH2013} Z. Ma, D. Callaway, and I. Hiskens, “Decentralized charging control of large populations of plug-in electric vehicles," IEEE Trans. Con. Sys. Tech., vol. 21, no. 1, pp. 67-78, 2013.

  \bibitem{basar2015} J. Moon and T. Basar. Discrete-time LQG mean field games with unreliable communication. 53rd IEEE Conference on Decision and Control, (2015), 2697-2702.

  \bibitem{basar2018} J. Moon and T. Basar. Risk-sensitive mean field games via the stochastic maximum principle. Dynamic Games and Applications,   (2018), 1-26.

    \bibitem{2}
      Y. Hu, J. Huang, and T. Nie. "Linear-quadratic-gaussian mixed mean-field games with heterogeneous input constraints." SIAM J. Control Optim., vol. 56, no. 4, pp. 2835–2877, 2018.


  \bibitem{NC2013} M. Nourian and P. E. Caines, “$\varepsilon$-Nash Mean Field Game Theory for Nonlinear Stochastic Dynamical Systems with Major and Minor Agents,” SIAM J. Control Optim., vol. 51, no. 4, pp. 3302–3331, Jan. 2013.

  \bibitem{NCMH2013}M. Nourian, P. E. Caines, R. P. Malhame and M. Huang, "Nash, Social and Centralized Solutions to Consensus Problems via Mean Field Control Theory," in IEEE Transactions on Automatic Control, vol. 58, no. 3, pp. 639-653, March 2013, doi: 10.1109/TAC.2012.2215399.





  \bibitem{QHX2020} Z. Qiu, J. Huang, and T. Xie, “Linear Quadratic Gaussian Mean-Field Controls of Social Optima,” arXiv preprint arXiv:2005.06792, 2019.

  \bibitem{RV1999} D. Ray and  R. Vohra, “A theory of endogenous coalition structures,” Games and economic behavior, vol. 26, no. 2, pp. 286-336, 1999.

  \bibitem{SLY16} J. Sun, X. Li, and J. Yong, “Open-Loop and Closed-Loop Solvabilities for Stochastic Linear Quadratic Optimal Control Problems,” SIAM Journal on Control and Optimization, vol. 54, no. 5. pp. 2274–2308, 2016.



  \bibitem{WH2017} B.-C. Wang and J. Huang, “Social optima in robust mean field LQG control,” 2017 11th Asian Control Conference (ASCC), 2017.


  \bibitem{YZ1999} J. Yong and X. Y. Zhou, Stochastic Controls: Hamiltonian Systems and HJB Equations. Springer Science \& Business Media, 1999.

\end{thebibliography}
\end{document}